\documentclass[11pt]{elsarticle}

\usepackage[usenames,dvipsnames]{xcolor}
\usepackage{amsmath,amsfonts,booktabs}
\usepackage{tikz,pgfplots}
\usepackage{graphicx,subfigure}
\usepackage{blkarray}
\usepackage{placeins}
\usepackage{multirow}
\usepackage{enumerate}
\usepackage{todonotes,setspace}
\usepackage{textcomp,gensymb,url}
\usepackage[left=1in, right=1.5in, top=1in, bottom=1in]{geometry}

\usetikzlibrary{external} 
\makeatletter
\def\ps@pprintTitle{%
 \let\@oddhead\@empty
 \let\@evenhead\@empty
 \def\@oddfoot{}%
 \let\@evenfoot\@oddfoot}
\makeatother

\pgfplotsset{mystyle/.style={mark=none}}

\begin{document}

\begin{frontmatter}
\title{A Fast and Memory Efficient Sparse Solver with Applications to Finite-Element Matrices} 


\author[Amir]{AmirHossein Aminfar\corref{cor1}\fnref{label1}}
\ead{aminfar@stanford.edu}
\author[Eric]{Eric Darve\fnref{label1}}

\fntext[label1]{Mechanical Engineering Department, Stanford University}
\cortext[cor1]{Corresponding author. +1 650-644-7624}
\address[Amir]{496 Lomita Mall, Room 104, Stanford, CA, 94305}
\address[Eric]{496 Lomita Mall, Room 209, Stanford, CA, 94305}


\begin{abstract}
In this article, we introduce a fast and memory efficient solver for sparse matrices arising from the finite element discretization of elliptic partial differential equations (PDEs). We use a fast direct (but approximate) multifrontal solver as a preconditioner, and use an iterative solver to achieve a desired accuracy. This approach combines the advantages of direct and iterative schemes to arrive at a fast, robust and accurate solver. We will show that this solver is faster ($\sim$ 2x) and more memory efficient ($\sim$ 2--3x) than a conventional direct multifrontal solver. Furthermore, we will demonstrate that the solver is both a faster and more effective preconditioner than other preconditioners such as the incomplete LU preconditioner. Specific speed-ups depend on the matrix size and improve as the size of the matrix increases. The solver can be applied to both structured and unstructured meshes in a similar manner. We build on our previous work and utilize the fact that dense frontal and update matrices, in the multifrontal algorithm, can be represented as hierarchically off-diagonal low-rank (HODLR) matrices. Using this idea, we replace all large dense matrix operations in the multifrontal elimination process with $O(N)$ HODLR operations to arrive at a faster and more memory efficient solver.
\end{abstract}

\end{frontmatter}


\section{Introduction}
\label {sec:intro}

In many engineering applications, we are interested in solving a set of linear equations:
\[
	Ax = b
\]
where $A$ is a symmetric positive definite stiffness matrix arising from a finite element discretization of an elliptic PDE, and $b$ is a forcing vector, associated with the inhomogeneity in the PDE. Iterative methods are widely popular in solving such equations. However, the main difficulty with these methods is that they require a preconditioner and convergence is not guaranteed. Direct methods on the other hand are very robust but are generally slower and more memory demanding. In this article, we present an accelerated multifrontal solver that we use as a preconditioner to a generalized minimal residual (GMRES~\cite{saad1986gmres}) iterative scheme to achieve a desired accuracy. This approach combines the robustness of direct solvers with the speed of iterative solvers to arrive at a fast overall solver for sparse finite element matrices. 

Accelerating the multifrontal direct solve algorithm has been the subject of many recent research articles~\cite{xia2009superfast,MFGeneralMesh,GeneralizedMF,randomizedMF,BLRMF}. For a detailed summary and overview of such algorithms see~\cite{BlackBox_HODLR}. The general idea behind most of these methods is approximating dense frontal matrices arising in the multifrontal elimination process with an off-diagonal low-rank matrix structure. The off-diagonal low-rank property leads to more efficient factorization and storage compared to dense BLAS3 operations if the rank is sufficiently small. The methods described in~\cite{xia2009superfast,MFGeneralMesh,GeneralizedMF,randomizedMF} approximate the frontal matrix with a hierarchically semiseparable (HSS) matrix, while \cite{BLRMF} approximates the frontal matrix with a block low-rank (BLR) matrix.

In this article, we accelerate the multifrontal algorithm by approximating dense frontal matrices as hierarchically off-diagonal low-rank (HODLR) matrices. Compared to HSS structures which have been widely used in approximating dense frontal matrices, HODLR matrices are much simpler as they lack the nested off-diagonal basis. For 3D PDEs, we find that the rank used to approximate the off-diagonal blocks increases with the size of the block with $r \approx O(\sqrt{n})$, where $r$ is the rank and $n$ the size of the block. This results in a geometric increase of the rank with the HODLR level. As a result of this increase, as we demonstrated in~\cite{BlackBox_HODLR}, the factorization cost is the same for both HODLR and HSS structures, namely $\mathcal{O}(r^2n)$, where $r$ is the rank at the top of the tree. This is despite the fact that HSS uses a more data-sparse format. The reason why the difference in the basis does not affect the asymptotic cost is because the cost is dominated by the computation at the root of the HODLR tree, for the largest block.

In addition, HODLR is advantageous compared to HSS during the low-rank approximation phase, since it does not need to produce a nested basis, which simplifies many steps in the algorithm. Hence, in most practical applications, HSS may not have a clear advantage over HODLR. 

Furthermore, we will demonstrate that the combination of HODLR and the boundary distance low-rank approximation method (BDLR)~\cite{BlackBox_HODLR} leads to a very fast and simple extend-add algorithm, which results in an overall fast multifrontal solver.

At the time of writing this article, only Xia~\cite{randomizedMF} has demonstrated a fast and memory efficient multifrontal solver for general sparse matrices, with an asymptotic cost in $O(N^{4/3} \log N)$ where $N$ is the size of the {\bf sparse} matrix. In contrast, the method in this paper leads to an overall cost of $O(N^{4/3})$. This cost may be compared with an LU factorization with nested dissection, with cost $O(N^2)$ in 3D.

In this article, we introduce a fast multifrontal solver that is much simpler compared to~\cite{randomizedMF}, and demonstrate its performance for large and complicated test cases. The method is shown to be advantageous compared to traditional preconditioners like ILU.

\section{Review of Important Concepts}

We now review two concepts that are central to the fast sparse solver algorithm. Namely, hierarchical off-diagonal low-rank (HODLR) matrices and the boundary distance low-rank approximation method (BDLR).

\subsection{Hierarchically Off-Diagonal Low-Rank (HODLR) Matrices}

Hierarchical matrices are data sparse representation of a certain class of dense matrices. This representation relies on the fact that these matrices can be sub-divided into a hierarchy of smaller block matrices, and certain sub-blocks can be efficiently represented as a low-rank matrix. We refer the readers to~\cite{hackbusch1999sparse, hackbusch2000sparse, grasedyck2003construction, hackbusch2002data, borm2003hierarchical,  ULV, chandrasekaran2006fast1,BlackBox_HODLR} for more details. Ambikasaran at al.~\cite{ambikasaran2013thesis} provides a detailed description of these different hierarchical structures. In this article, we use the simplest hierarchical structure, namely the hierarchically off-diagonal low-rank matrix (HODLR), to approximate the dense frontal matrices that arise during the sparse elimination process. As shown in~\cite{BlackBox_HODLR}, the HODLR structure reduces the dense factorization and storage cost from $\mathcal{O}(n^3)$ and $\mathcal{O}(n^2)$ to $\mathcal{O}(r^2n)$ and $\mathcal{O}(rn)$ respectively, where $n$ is the size of the dense matrix and $r$ is the off-diagonal rank.

An HODLR matrix has low-rank off-diagonal blocks at multiple levels. As described in~\cite{SivaFDS}, a 2-level HODLR matrix, $K\in \mathbb{R}^{n\times n}$, can be written as shown in Eq.~\eqref{eq:HODLR2}:
\begin{align}
	K & =
	\begin{bmatrix}
		K_1^{(1)}&U_1^{(1)} (V_{1,2}^{(1)})^T \\
		U_2^{(1)} (V_{2,1}^{(1)})^T&K_2^{(1)}
	\end{bmatrix} \notag \\
	& =
	\begin{bmatrix}
		\begin{bmatrix}
			K_1^{(2)}&U_1^{(2)} (V_{1,2}^{(2)})^T \\
			U_2^{(2)} (V_{2,1}^{(2)})^T&K_2^{(2)}
		\end{bmatrix}&
		         U_1^{(1)} (V_{1,2}^{(1)})^T \\
			U_2^{(1)} (V_{2,1}^{(1)})^T&
		\begin{bmatrix}
			K_3^{(2)}& (U_3^{(2)})^T (V_{3,4}^{(2)})^T \\
			U_4^{(2)} (V_{4,3}^{(1)})^T&K_4^{(2)}
		\end{bmatrix}
	\end{bmatrix}
\label{eq:HODLR2}	
\end{align}
where for a $p$-level HODLR matrix, $K_i^{(p)} \in \mathbb{R}^{n/2^p\times n/2^p}$, $U_{2i-1}^{(p)}$, $U_{2i}^{(p)}$, $V_{2i-1,2i}^{(p)}$, $V_{2i,2i-1}^{(p)} \in \mathbb{R}^{n/2^p\times r}$ and $r\ll n$. Further nested compression of the off-diagonal blocks will lead to an HSS structure~\cite{SivaFDS}.

\subsection{Boundary Distance Low-Rank Approximation Method (BDLR)}
\label{sec:BDLR}
In order to take advantage of the off-diagonal low-rank property, we need a fast and robust low-rank approximate method. More precisely, we need a low rank approximation method that has the following properties:
\begin{itemize}
	\item We want our method to be applicable to general sparse matrices. Hence, we need a low-rank approximation scheme that is purely algebraic (black-box). That is, we can not use analytical low-rank approximation methods like Chebyshev, multipole expansion, analytical interpolation, etc.
	\item In order to obtain speedup compared to conventional multifrontal solvers, we need a fast low-rank approximation scheme that has a computational cost of $\mathcal{O}(rn)$ where $n$ and $r$ are the size and rank of a dense low-rank matrix respectively. Hence, we cannot use traditional low-rank approximation methods like SVD, rank revealing LU or rank revealing QR as they have a computational cost of $\mathcal{O}(n^3)$, $\mathcal{O}(n^2)$ and $\mathcal{O}(n^2)$ respectively. 
	\item We need a robust and efficient low-rank approximation method that is applicable to a wide variety of problems.
\end{itemize}
One possible option is to use randomized algorithms~\cite{RndSummary,Rnd1,Rnd2,Rnd3} similar to Xia~\cite{randomizedMF}. However, such algorithms require the implementation of a fast matrix-vector product. For our purpose, as we demonstrated in~\cite{BlackBox_HODLR}, the boundary distance low-rank approximation method (BDLR) is a fast and robust scheme that results in very fast solvers for both structured and unstructured meshes. 

\begin{figure}[htbp]
	\centering
	\subfigure[Full Matrix Representation]{
	\includegraphics[scale=1]{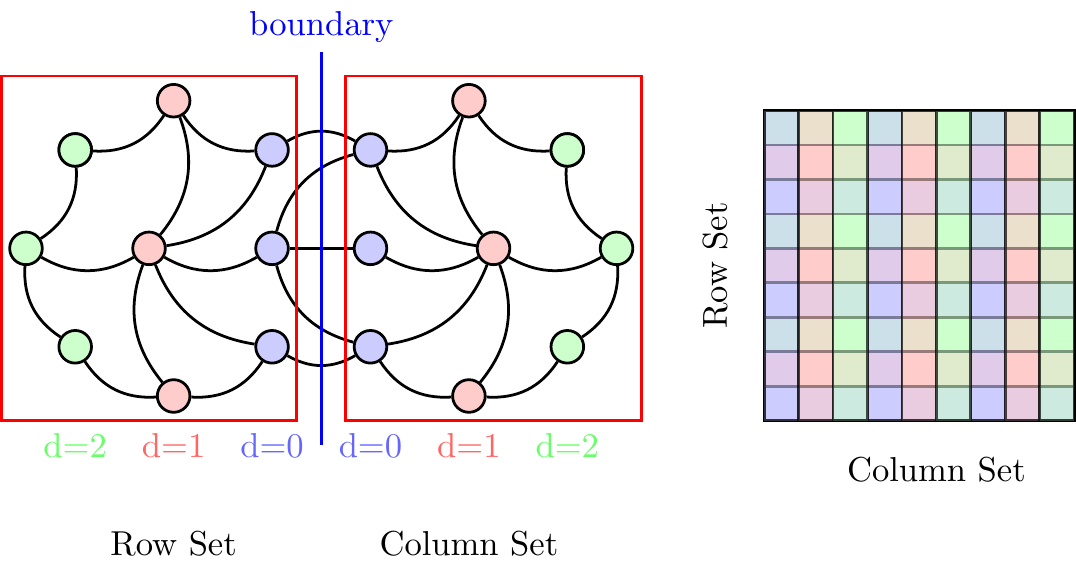}
	\label{fig:BDLR_Full}
	}
	\subfigure[Low-Rank Matrix Representation]{
	\includegraphics[scale=1]{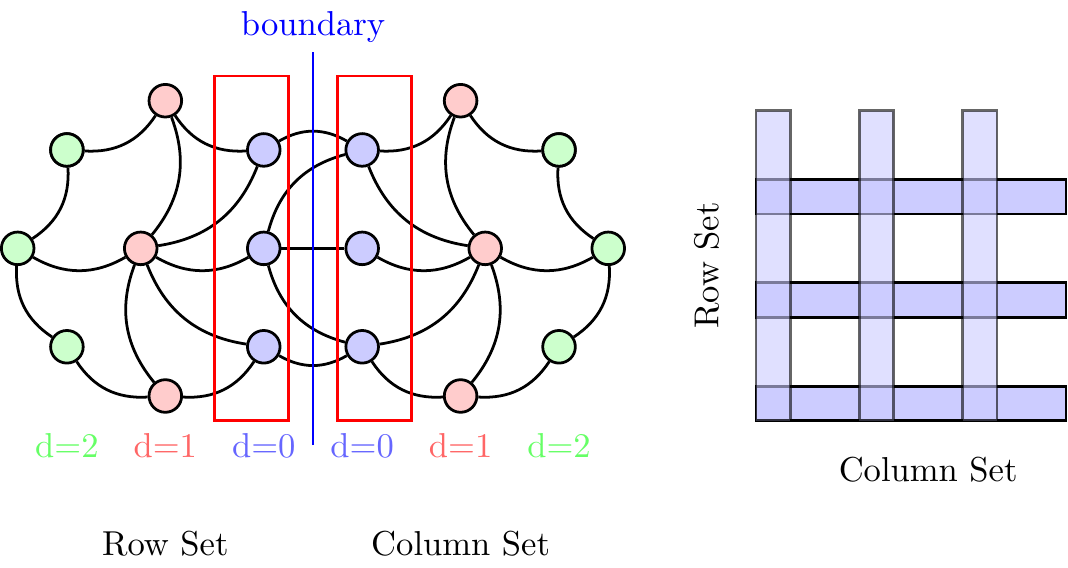}
	\label{fig:BDLR_LR}
	}
	\caption{Classification of vertices based on distance from the other set.}
\end{figure}

BDLR is a pseudoskeleton~\cite{pseudoSkeleton} like low-rank approximation scheme that picks rows and columns based on the corresponding interaction graph of a dense matrix, which in the case of frontal matrices, is the graph corresponding to the sparse separator. That is, for an off-diagonal block in the frontal matrix, it chooses a subset of rows and columns  based on the corresponding separator graph. The criteria for choosing these rows and columns is based on the location of their respective nodes in the sparse separator graph. Figure~\ref{fig:BDLR_Full} shows an example of an interaction graph corresponding to the interaction of a set of row and column indices in an off-diagonal block of a sample frontal matrix. Figure~\ref{fig:BDLR_LR} shows that the BDLR method chooses row indices and column indices corresponding to nodes that are closer to the boundary (blue line). 

Let $R$ and $C$ be the matrices containing all the selected rows and columns. In other words:
\[
R = B(I,:), \qquad
C = B(:,J)
\]
where $B$ is the off-diagonal low-rank matrix and $I$ and $J$ are the set of row and columns indices chosen by the BDLR algorithm respectively. Defining $\widehat{B} = B(I,J)$, we perform a full pivoting LU factorization :
\begin{equation}
	\label{eq:Bhat}
	{\widehat{B}} = P^{-1}LUQ^{-1}
\end{equation}
where $P$ and $Q$ are permutation matrices. Let $r$ be the chosen rank for $\widehat{B}$. Define $\widetilde{R}$ and $\widetilde{C}$ as:
\begin{align*}
  \widetilde{C}& =  (CQ)(:,1:r) \; (U(1:r,1:r))^{-1}\\
  \widetilde{R}& =  (L(1:r,1:r))^{-1} \; (PR)(1:r,:)
\end{align*}
We then have:
\begin{equation*}
	B\approx \widetilde{C} \; \widetilde{R}
\end{equation*}
$(U(1:r,1:r))^{-1}$ and $(L(1:r,1:r))^{-1}$ correspond to lower-triangular solves. The inverse matrices are not explicitly computed. The approximation rank $r$ is chosen based on the desired final accuracy such that $|u_{r+1,r+1} / u_{11}| < \epsilon$ where $u_{r+1,r+1}$ and $u_{11}$ correspond to $(r+1)$th and the first pivots respectively and $\epsilon$ is the desired accuracy. 

The final rank $r$ may be significantly smaller than the number of originally selected rows and columns. This higher compression results in higher efficiency both in terms of memory and runtime.

\section {An Iterative Solver with Direct Solver Preconditioning}
\label{sec:directIterative}

In this article, we investigate using an accelerated multifrontal sparse direct solver as a preconditioner to the generalized minimal residual (GMRES)~\cite{saad1986gmres} method. In this case, we use a relatively low accuracy for the direct solver. We will show that this approach is much faster and more memory efficient than a conventional multifrontal sparse solver. We should also mention that this preconditioning method can be applied to any iterative solvers (conjugate gradient (CG)~\cite{hestenes1952methods}, etc.). 

\section {A Fast Multifrontal Solver}
\label{sec:directSolver}

\subsection{Overview of a Conventional Multifrontal Algorithm}
\label{sec:conventional}

We do not intend to give a detailed explanation of the multifrontal solve process in this article. We refer the reader to the available literature (see for example~\cite{MFReview}) for an in-depth explanation of the algorithm.

In the multifrontal method, the unknowns are eliminated following the ordering imposed by the elimination tree. That is, each node of the elimination tree corresponds to a set of unknowns, and these unknowns cannot be eliminated until all the unknowns corresponding to the children of this node are eliminated.

The multifrontal algorithm is an algorithm to calculate the Cholesky or LU factorization of a sparse matrix~\cite{MFReview}, with special optimizations that take advantage of the sparsity. Moreover (and this is specific to a multifrontal elimination), during the elimination, information is propagated only from a child node to its parent (in the so-called elimination tree~\cite{duff1986direct}). This is what distinguishes for example a multifrontal elimination from a supernodal elimination.

We note that in this paper we describe our method in the context of a multifrontal elimination; however, the same method can be applied to a supernodal elimination. No fundamental change is required to our algorithm.

\subsubsection{Factorization}
Consider now a node $p$ in the elimination tree. Let $I_p$ be the set of indices of unknowns associated with node $p$:
\begin{equation*}
 I_p = \{i^{(p)}_{1}, \ldots, i^{(p)}_{n_p} \}
\end{equation*}
where $n_p$ is the number of unknowns corresponding to node $p$, and $i^{(p)}_j$ is the global index of the $j$th unknown associated with node $p$. We denote a specific child node of $p$ as $c_k$ ($c_k \in \mathcal{C}_p$, $k \in \{1, \ldots,n^c_p \}$ where $\mathcal{C}_p$ is the set of all children and $n^c_p$ is the number of children of node $p$).

Define the set $S_p$ as the set of unknowns $j > i^{(p)}_{n_p}$  that are connected to any of the unknowns in $I_p$, in the graph of $A$. More precisely:
\begin{equation}
\label{eq:couplingSet}
S_p = \{ j \, | \, \exists i \in I_p, j > i^{(p)}_{n_p}, a_{ij} \ne 0 \}
\end{equation}
where $a_{ij}$ is the entry at the $i$th row and $j$th column of the original sparse matrix $A$. Describing the details of a multifrontal elimination requires the definition of the matrix $U_{c_k}$, which we call the update matrix corresponding to the $k$th child of node $p$ by recurrence. The set of indices corresponding to unknowns associated with $U_{c_k}$ (update matrix of children nodes) is denoted $I^U_{c_k}$. 

If $p$ is a leaf node in the elimination tree, the matrix $U_{c_k}$ is not defined and hence, $I^U_{c_k}=\emptyset$. We define the set of frontal indices $I_p^f$ as follows:
\begin{equation*}
I_p^f = S_p \cup \{\cup_{k = 1}^{n^c_p} I^U_{c_k} \} \setminus I_p
\end{equation*}
We now define the matrix $\bar{F}_p$ as the sub-matrix of $A$ associated with $I_p \cup I_p^f$.
\begin{equation}
	\label{eq:fbar}
	\includegraphics{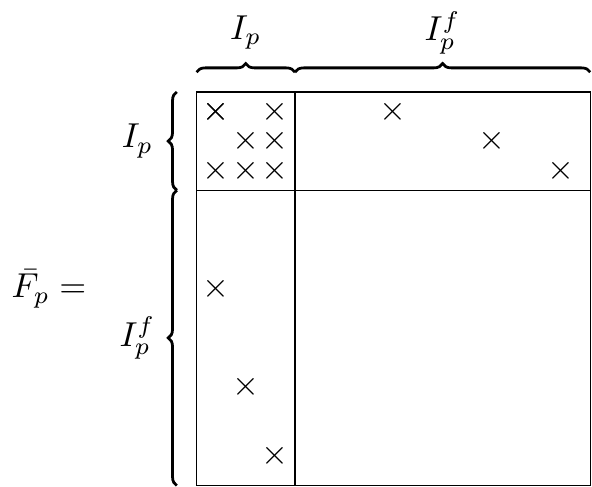}
\end{equation}
The symbol $\times$ schematically denotes a nonzero entry in the matrix. In $\bar{F}_p$, we set the entries for the block $I_p^f \times I_p^f$ to 0.

The frontal matrix for node $p$, $F_p$, is defined as follows:
\begin{equation*}
	F_p = \bar{F}_p \oplus U_{c_1} \oplus \cdots \oplus U_{c_{n^c_p}}	
\end{equation*}
The symbol $\oplus$ denotes the extend-add operation and $U_{c_k}$ denotes the update matrix corresponding to the $k$th child of node $p$. The extend-add is basically an addition. The ``extend'' part corresponds to the fact that there is a size mismatch between a $U_{c_k}$ and $F_p$, and an index mapping from $U_{c_k}$ to $F_p$ must be used. In the special case where node $p$ is a leaf node in the elimination tree, $F_p = \bar{F}_p$.

Note that, after the extend add operations, the frontal matrix $F_p$ is (nearly) a fully dense matrix. We then divide $F_p$ into four sub-blocks:
\begin{equation*}
	\includegraphics{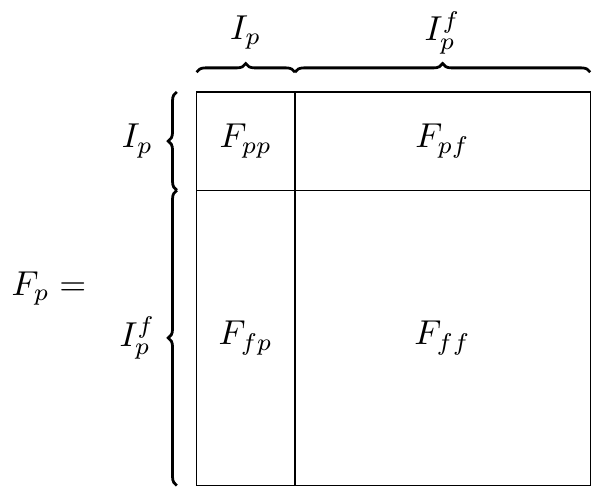}
\end{equation*}
Factorizing $F_{pp}$, we are left with the Schur complement. This is by definition the update matrix associated with node $p$ which will be used in the extend-add operation of its parent:
\begin{equation}
\label{eq:update}
U_p \overset{\text{def}}{=} F_{ff} - F_{fp}F_{pp}^{-1}F_{pf}
\end{equation}
Repeating the operations described in Eqns.~\eqref{eq:couplingSet} to~\eqref{eq:update} for all nodes in the elimination tree starting from the leaf nodes and going up to the root node, constitutes the factorization phase of the conventional multifrontal algorithm.

\subsubsection{Solve}
\label{sec:convSolve}
The solve phase constitutes of an upward pass (L solve) and a downward pass (U solve) in the elimination tree. In the upward (downward) pass, we traverse the elimination tree upward (downward) from leaves to root (root to leaves) and traverse the right hand side vector $b$ downwards (upwards). Hence, the upward and downward passes correspond to the L and U solve phases in a conventional LU solver respectively.

In the upward pass (L solve) phase, we first construct the upward pass solution matrix $b_u$ which is initially equal to the right hand side $b$. Then, moving upward in the elimination tree, we construct the upward solution $b_{u_p}$ for each node $p$, which is basically the elements of the upward pass solution $b_u$ corresponding to the unknowns in $I_p$ and $I_p^f$.
\begin{equation}
	\label{eq:L_Soln}
	\includegraphics{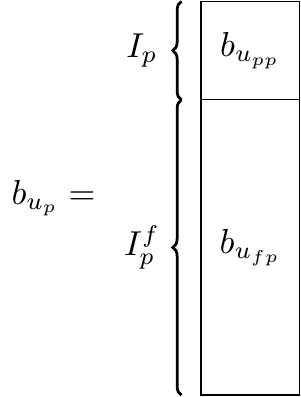}
\end{equation}
Now, update the upward pass solution using:
\begin{equation}
	\label{eq:upward}
	b_{u_{fp}} = b_{u_{fp}} - F_{fp}F_{pp}^{-1}b_{u_{pp}}
\end{equation}

After completing the upward pass (L solve), we must perform a downward pass (U solve) to arrive at our final solution. The final solution $x$ is initially an empty vector. We traverse the elimination tree from root to leaves (downward). For each node $p$, we construct the final solution vector (Eq.~\eqref{eq:U_Soln}). The corresponding solution for each node can be calculated as follows:
\begin{gather}
	\label{eq:U_Soln}
	\includegraphics{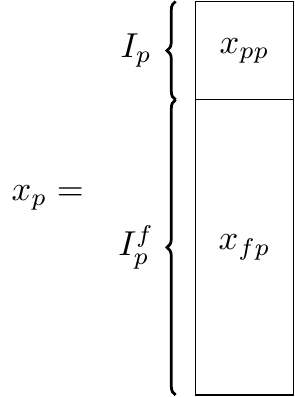} \\
	\label{eq:downward}
	x_{pp} =  F_{pp}^{-1}(b_{u_{pp}} - F_{pf}x_{fp})
\end{gather}
Note that since we're traversing the elimination tree downward (traversing $b_{u_p}$ upward), $x_{fp}$ has already been calculated by the time we reach $p$.

\subsection{HODLR Accelerated Multifrontal Algorithm}
Looking at the procedure described in Section~\ref{sec:conventional}, one can observe that dense BLAS3 operations like the one described by Eq.~\eqref{eq:update}, which involves both a factorization and an outer product update, can become time and memory consuming as the front size increases. In order to accelerate the multifrontal elimination process, we replace large dense matrices with HODLR structures.

\subsubsection{Accelerated Factorization}
In the factorization phase, we want to represent the frontal and update matrices for each node $p$ as HODLR matrices. In order to be able to construct an HODLR structure, we need to utilize a suitable low-rank approximation method. Our previous results~\cite{BlackBox_HODLR} show that the boundary distance low-rank approximation scheme (BDLR) is a suitable algorithm for our purposes. Furthermore, as we will show, a priori knowledge of rows and columns with BDLR leads to a very fast extend-add operation.

To construct the HODLR representation of the frontal and update matrices of $p$, we first assemble $\bar{F}_p$ as described by Eq.~\eqref{eq:fbar}. As described in~\cite{BlackBox_HODLR}, the BDLR algorithm requires an interaction graph that describes the interaction between the rows and columns of the matrix, which, in this case, is the graph that is constructed from the submatrix of the original matrix $A$ corresponding to the interaction of rows and columns with indices in the set $I_p \cup I_p^f$. 

Using the interaction graph for $I_p \cup I_p^f$, we create the HODLR representation of $F_p^{HODLR}$ using the BDLR algorithm:
\begin{equation}
	\label{eq:frontalMatrix_HODLR}
	\includegraphics[scale=1.1]{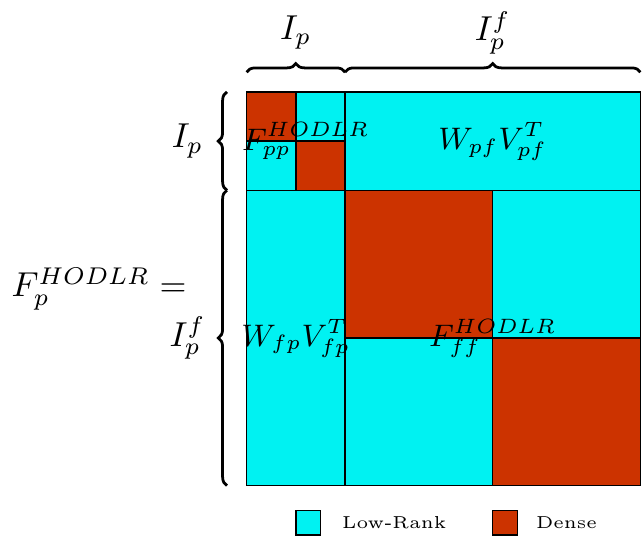}
\end{equation}
Using the extend-add notation $\oplus$, $F_p^{HODLR}$ is given by:
\begin{equation}
	\label{eq:HODLR_EA}
	F_p^{HODLR} = \bar{F}_p \oplus U_{c_1}^{HODLR} \oplus \cdots 
	\oplus U_{c_{n_c^p}}^{HODLR}	
\end{equation}
For simplicity, we've assumed that all the update matrices associated with node $p$ are HODLR matrices. In some cases, the update matrices might be small dense matrices, in which case the extend-add operations described for the child HODLR updates will become almost trivial for the dense child updates. Looking at Eq.~\eqref{eq:update}, we notice that every HODLR update matrix is composed of two components: an HODLR matrix and an outer product update.
\begin{equation}
\label{eq:update_HODLR}
U_{c_k}^{HODLR} = F_{{ff}_k}^{HODLR} - W_{{fp}_k}V_{{fp}_k}^{T}(F_{{pp}_k}^{HODLR})^{-1}W_{{pf}_k}V_{{pf}_k}^T=F_{{ff}_k}^{HODLR} - W_kV_k^T
\end{equation}
where the subscript $k$ denotes the update from the $k$th child of $p$. Since we only utilize $U_{c_k}^{HODLR}$ in the extend-add operation of its parent, we will save the two contributions for the extend-add operation. That is, Eq.~\eqref{eq:HODLR_EA} now becomes:
\begin{equation}
	\label{eq:HODLR_EA_Expand}
	F_p^{HODLR} = \bar{F}_p \oplus (F_{{ff}_1}^{HODLR} - W_1V_1^T) \oplus \cdots \oplus (F_{{ff}_{n_c^p}}^{HODLR} - W_{n_c^p}V_{n_c^p}^T)
\end{equation}
Equation~\eqref{eq:HODLR_EA_Expand} requires that we perform an extend-add operation into a target HODLR structure. 

Before going into the details of this operation, we should first emphasize the importance of this operation both in terms of computational cost and memory saving compared to a conventional extend-add operation. Consider the outer product $W_1V_1^T$ in Eq.~\eqref{eq:HODLR_EA_Expand}. Let $\hat{n}_p$ be the size of the matrix $F_p^{HODLR}$. In order to perform the extend-add operation, we must extend $W_1$ and $V_1$ to arrive at matrices $W_1^e$ and $W_2^e$, of size $\hat{n}_p\times r_p$. In the conventional algorithm, we had to perform the outer product $W_1^eV_1^{e^T}$, which has a computational cost of $\mathcal{O}(\hat{n}_p^2r_p)$ and a storage cost of $\mathcal{O}(\hat{n}_p^2)$. For a 3D mesh with $N$ degrees of freedom, $\widehat{n}_p$, corresponding to the root node in the elimination tree, grows as $\mathcal{O}(N^{2/3})$, and $r_p$ for this node roughly scales as $\mathcal{O}(N^{1/3})$. This will result in a computational cost of $\mathcal{O}(N^2)$ and a storage cost of $\mathcal{O}(N^{4/3})$. Hence, in practice, the extend-add operation dominates the computational cost of the conventional multifrontal algorithm. 

\begin{figure}[htbp]
	\centering
	\subfigure[Determine the required off-diagonal rows and columns using BDLR]{
		\includegraphics[scale=0.6]{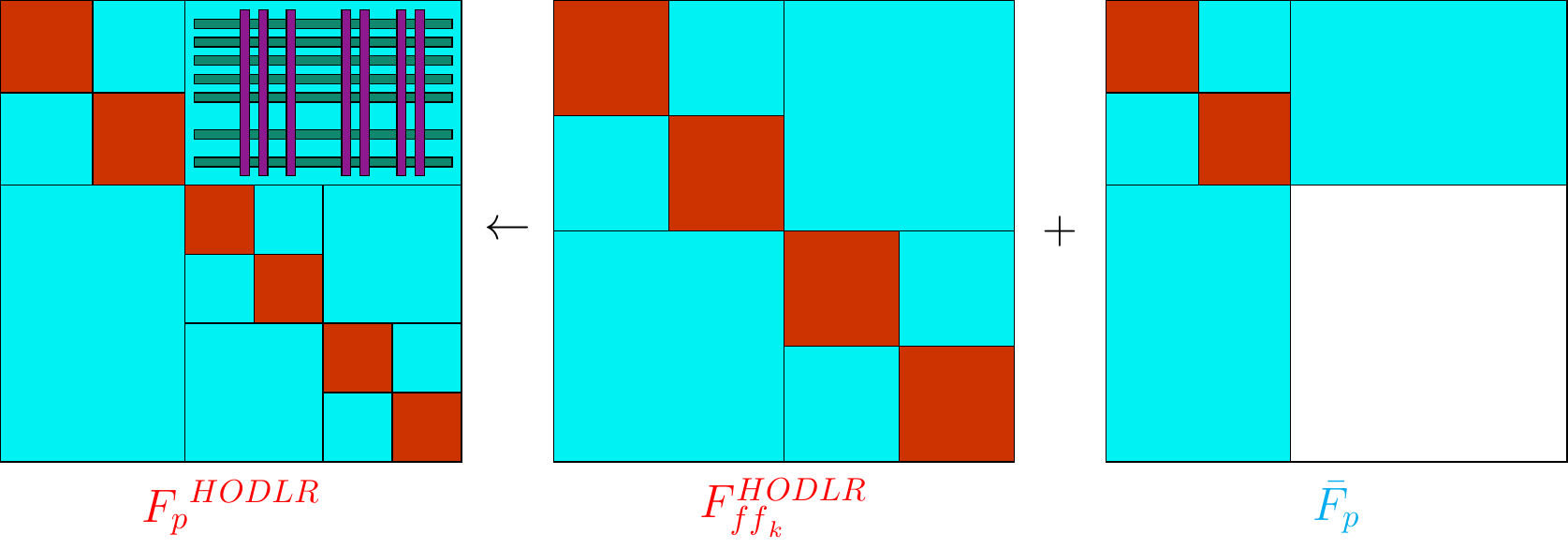}
	}
	\subfigure[Extract the identified rows and columns from $\bar{F_p}$]{
		\includegraphics[scale=0.6]{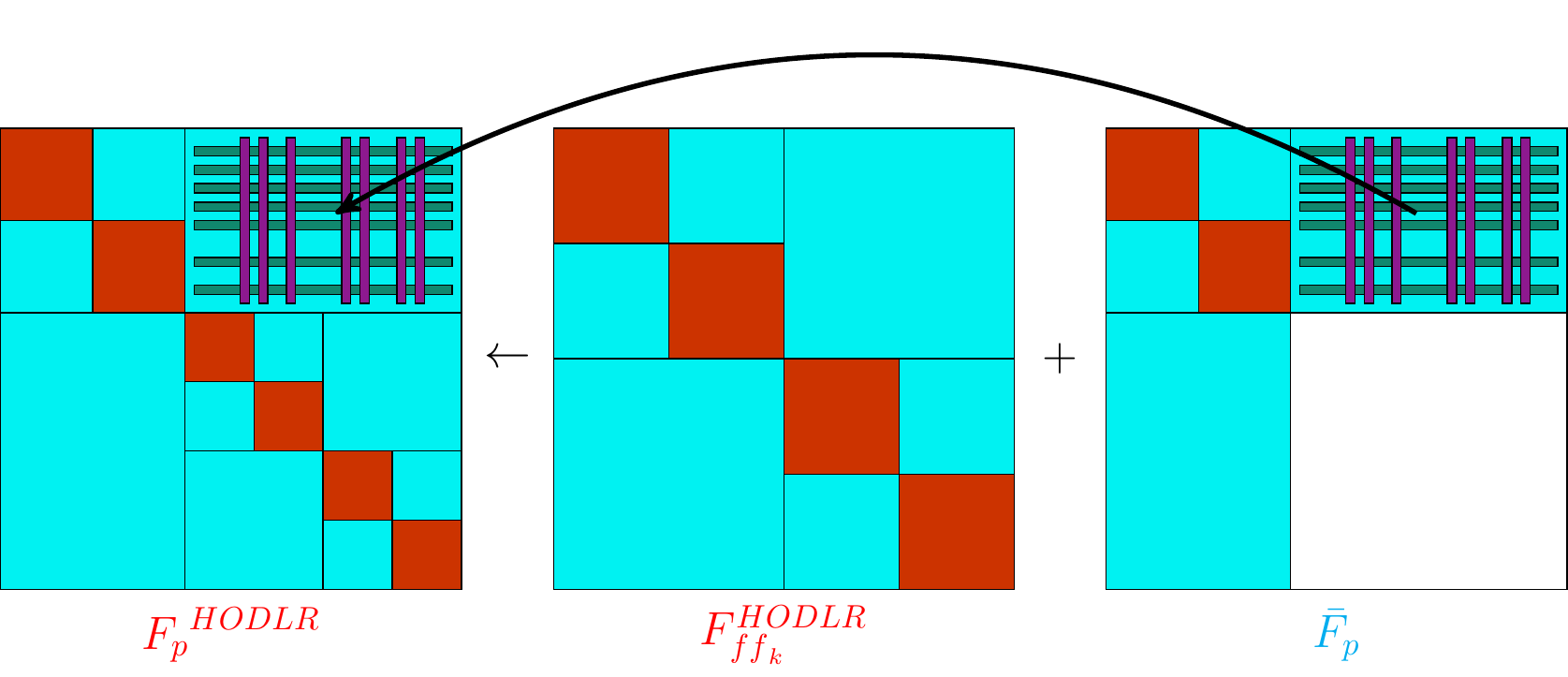}
	}
	\subfigure[Extract the identified rows and columns from $F_{{ff}_k}^{HODLR}$]{
		\includegraphics[scale=0.6]{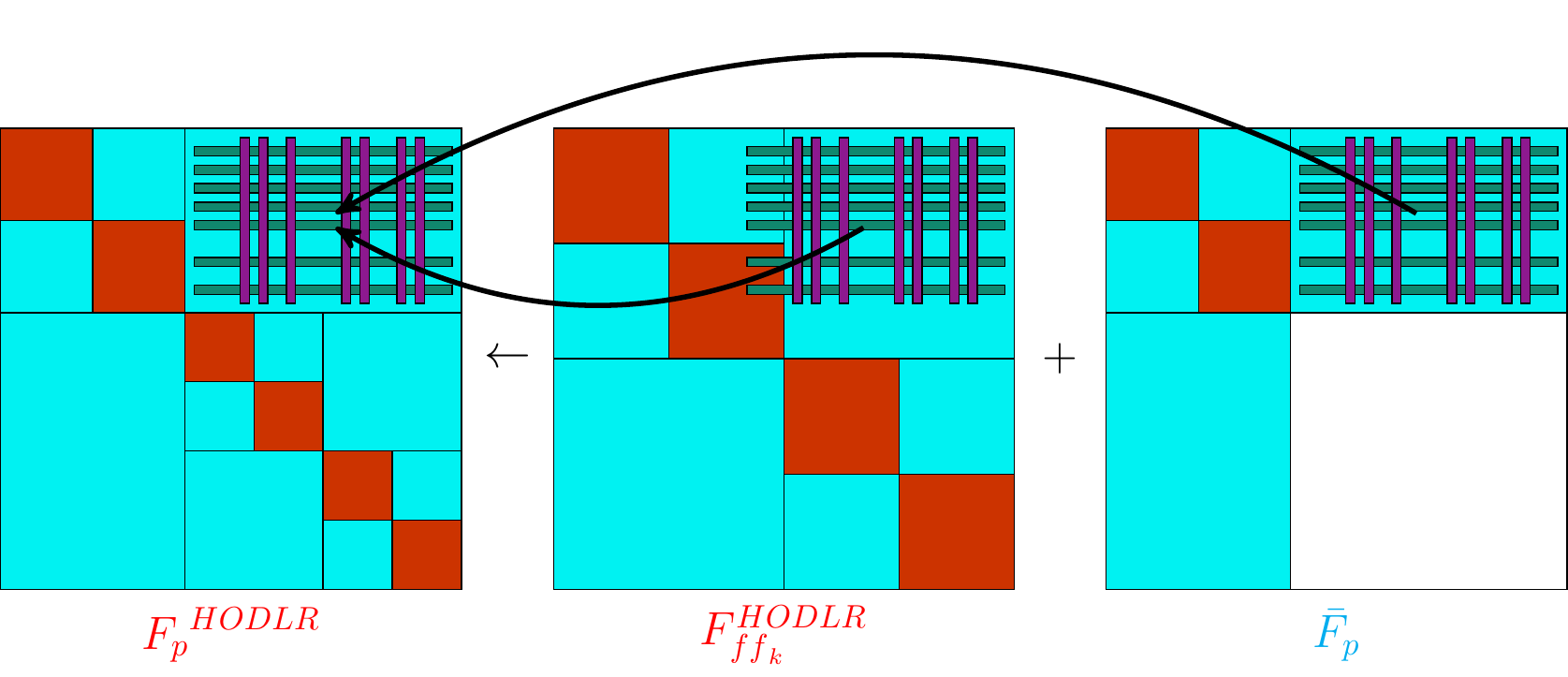}
	}
	\subfigure[Repeat the same procedure for all off-diagonal blocks of the target matrix]{
		\includegraphics[scale=0.6]{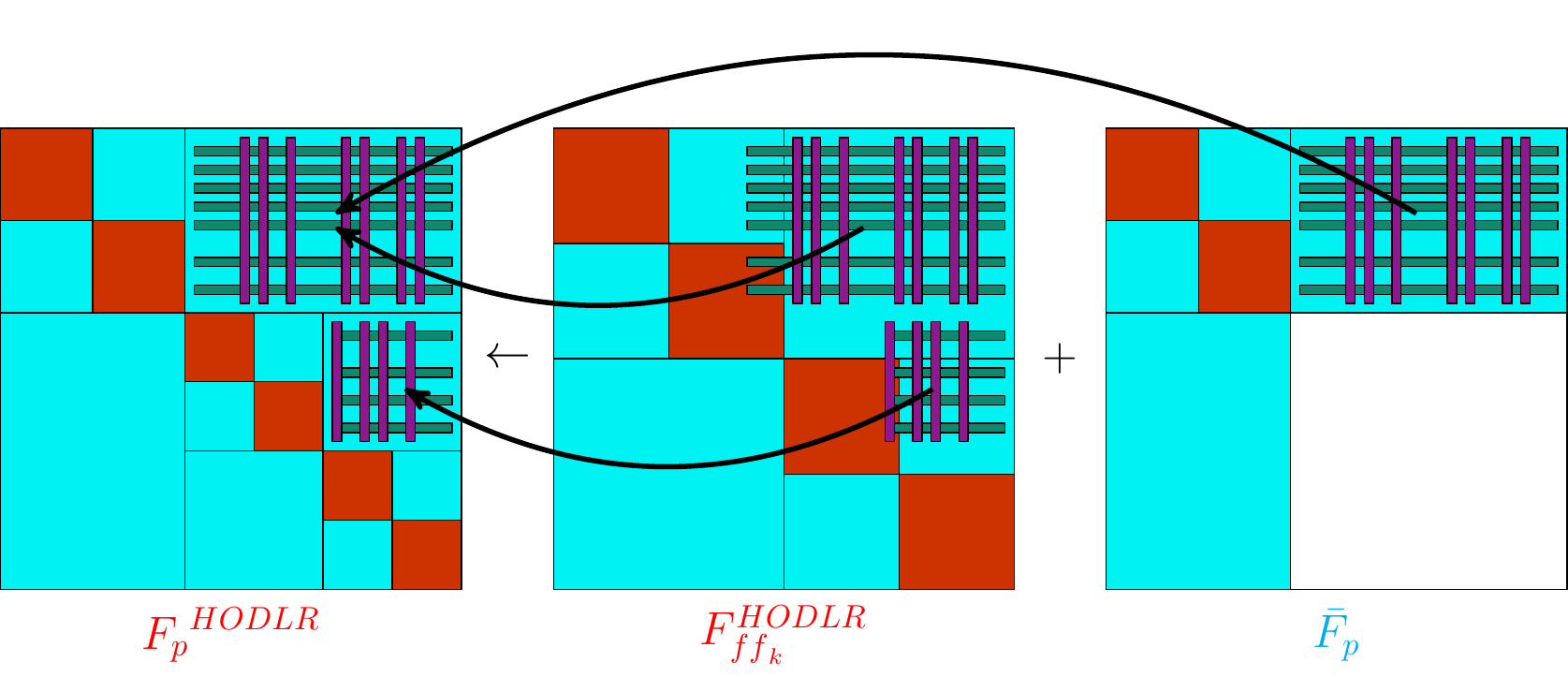}
	}
\caption{Fast HODLR $\leftarrow$ HODLR $+$ HODLR operation using the BDLR low-rank approximation algorithm. Red: Dense matrix block, Cyan: Low-rank matrix block, White: Block of zeros. \label{fig:HODLR_EA}}
\end{figure}

As shown in Eq.~\eqref{eq:HODLR_EA_Expand}, the extend-add process involves two different operations. The first operation is updating the frontal matrix ($\bar{F}_p$) with an HODLR structure ($F_{{ff}_k}^{HODLR}$) to arrive at a target HODLR structure ($F_p^{HODLR}$). What makes this operation difficult is the fact that $F_{{ff}_k}^{HODLR}$ and $F_p^{HODLR}$ typically have different structures. See Figure~\ref{fig:HODLR_EA} for an illustration. That is, the diagonal block sizes and the number of HODLR levels might differ between the two matrices. 

A key feature of the BDLR algorithm is that for each target HODLR structure ($F_p^{HODLR}$), we know a priori the rows and columns needed to construct the off-diagonal low-rank approximation. Hence, in order to perform the extend-add operation, we traverse the target HODLR structure and for each off-diagonal block, we extract the rows and columns determined by the BDLR algorithm from the child update HODLR matrices ($F_{{ff}_k}^{HODLR}$).

The second extend add-operation is adding a low-rank matrix ($W_kV_k^T$) to the frontal matrix ($\bar{F}_p$) to arrive at a target HODLR structure ($F_p^{HODLR}$). This process is very similar to the one described for adding $F_{{ff}_k}^{HODLR}$ to $\bar{F}_p$. The only difference is that instead of reconstructing rows and columns from an HODLR structure ($F_{{ff}_k}^{HODLR}$), we reconstruct the required rows and columns from the low-rank outer product ($W_kV_k^T$).

After reconstructing and extracting the rows and columns selected by BDLR and adding them to the corresponding rows and columns in the target structure, we perform the partial pivoting LU factorization. As described in Section~\ref{sec:BDLR}, we arrive at a final rank $r$ which is much smaller compared to the number of originally selected rows and columns. Given that the factorization of a hierarchical matrix of size $n$ scales as $\mathcal{O}(r^2n)$, a reduction in the rank has a significant effect on the resulting speedup.

Assuming the off-diagonal rank of the target HODLR structure corresponding to the node $p$ in the elimination tree is $r_p$, we need to extract $r_p$ rows and columns from the HODLR matrices and the outer products. Hence, we need to perform $\mathcal{O}(r_p^2\hat{n}_p)$ operations in order to construct $r_p$ rows from the outer product updates. This translates to a computational cost of $\mathcal{O}(N^{4/3})$  for the root node in the elimination tree ($r_p$ scales as $\mathcal{O}(N^{1/3})$) and is much more efficient compared to the $\mathcal{O}(N^2)$ scaling of the conventional extend-add algorithm. Moreover, we used no additional memory in order to perform the accelerated extend-add operation.

Now that we have constructed the HODLR representation of the frontal matrix $F_p^{HODLR}$, we factorize $F_{{pp}}^{HODLR}$ using an HODLR solver (see~\cite{BlackBox_HODLR} for example). Next, we store the update matrix $U_p^{HODLR}$ as an HODLR matrix and an outer product update:
\begin{equation}
\label{eq:update_HODLR_Parent}
U_p^{HODLR}\overset{\text{def}}{=} F_{{ff}}^{HODLR} - W_{{cp}}V_{{cp}}^{T}(F_{{ff}}^{HODLR})^{-1}W_{{pc}}V_{{pc}}^T=F_{{ff}}^{HODLR} - W_pV_p^T
\end{equation}

\subsubsection{Accelerated Solve}
The solve phase of the accelerated method is very similar to the solve phase of the conventional method described in Section~\ref{sec:convSolve}. The only difference is that $F_{pp}^{-1}$ is now replaced by $(F_{pp}^{HODLR})^{-1}$, which simply represents an HODLR solve instead of a conventional solve. Furthermore, the matrices $F_{pf}$ and $F_{fp}$ are now represented as low-rank products which results in a more efficient matrix-vector multiplications in Eqns.~\eqref{eq:upward} and~\eqref{eq:downward}.

\section {Application to Finite-Element Matrices}

In order to demonstrate the effectiveness of our method, we benchmark our solver for two classes of problems. We first apply our solver to a finite-element stiffness matrix that arises from a complicated 3D geometry. Next, we benchmark the performance of our solver for sparse matrices arising from the FETI method~\cite{FETI_DP1,FETI_DP2}.

\subsection{Elasticity Problem for a Cylinder Head Geometry}
\label{sec:cylinderResults}

\begin{figure}[htbp]
\centering
	\includegraphics[width=200pt,height=165pt]{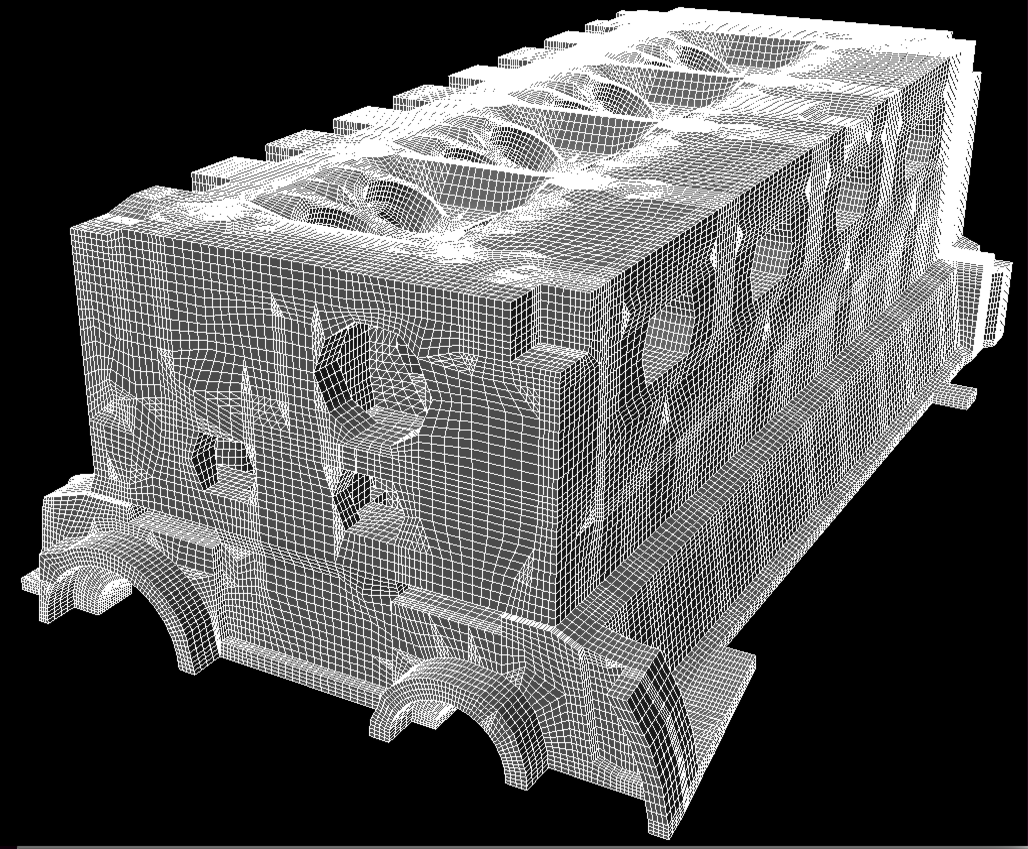}
\caption{A sample cylinder head mesh.}
\label{fig:cylinderMesh}
\end{figure}

We apply the iterative solver with the accelerated multifrontal preconditioner to a stiffness matrix corresponding to the finite-element discretization of the elasticity equation in a cylinder head geometry:
\begin{equation}
	(\lambda + \mu)\nabla(\nabla \cdot \boldsymbol{u})+\mu \nabla^2\boldsymbol{u}+\boldsymbol{F} = 0
	\label{eq:NavierCauchy}
\end{equation}
where $\boldsymbol{u}$ is the displacement vector and $\lambda$ and $\mu$ are Lam\'e parameters. The cylinder head mesh consists of a mixture of 8-node hexahedral, 6-node pentahedral and 4-node tetrahedral solid elements, and also 3-node shell elements. Figure~\ref{fig:cylinderMesh} shows a sample mesh for the cylinder head geometry.

\subsection{FETI-DP Solver for a 3D Elasticity Problem}
\label{sec:FETIResults}

FETI methods~\cite{FETI_DP1,FETI_DP2} are a family of domain decomposition algorithms with Lagrange multipliers that have been developed for the fast sequential and parallel iterative solution of large-scale systems of equations arising from the finite-element discretization of partial differential equations~\cite{FETI_DP1}. In this article, we investigate the solution of sparse matrices arising from a FETI-DP solver applied to the elasticity equation Eq.~\eqref{eq:NavierCauchy}.

We consider two classes of problems within the FETI-DP framework. The first class of matrices, called local matrices, corresponds to solving the problem on a subdomain of the original mesh. The other class of matrices is called coarse problem matrices and corresponds to the corner DOFs of all the subdomains.

\begin{figure}[htbp]
	\centering
	\subfigure[structured cube]{
		\includegraphics[width=190pt]{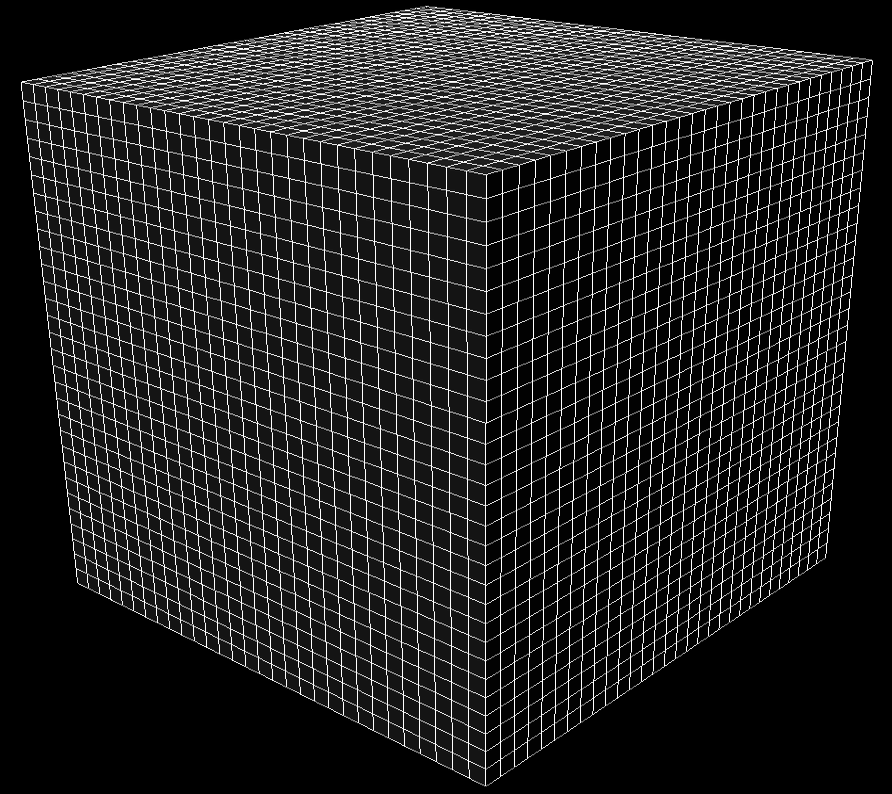}
	\label{fig:structuredMesh}
	}
	\subfigure[unstructured cube]{
		\includegraphics[width=190pt]{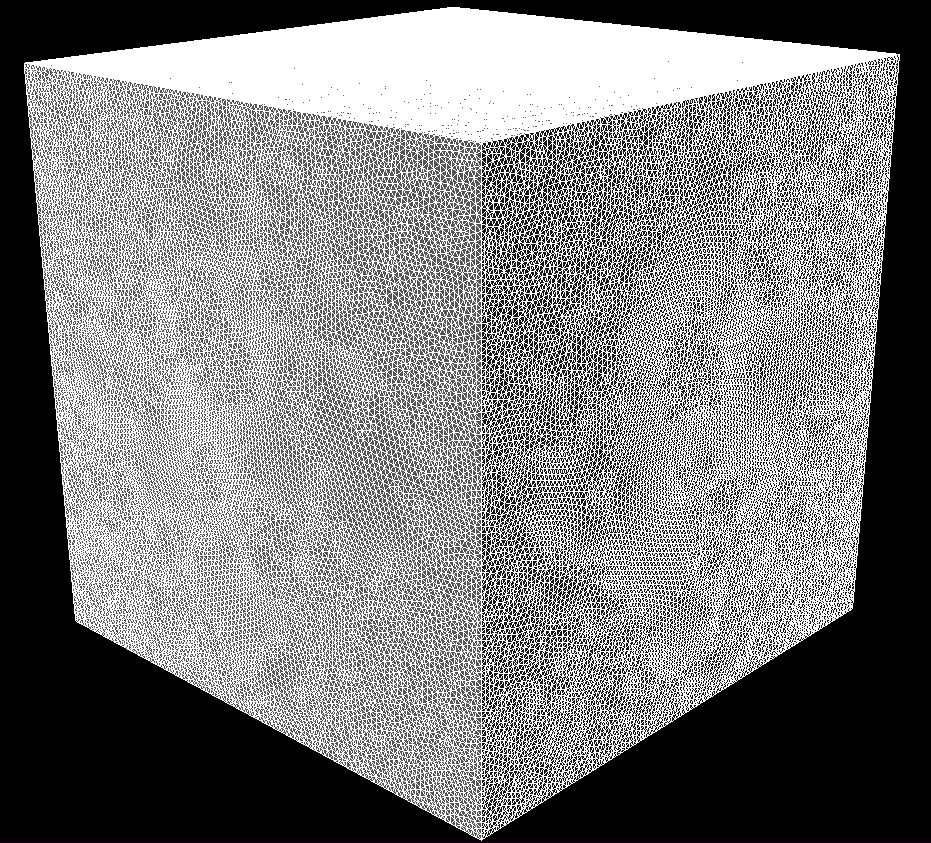}
	\label{fig:unstructuredMesh}
	}
\caption{Sample structured cube, unstructured cube and vehicle meshes. These meshes correspond to solving the local FETI-DP problem corresponding to the finite-element discretization of the elasticity equation Eq.~\eqref{eq:NavierCauchy}.}
\end{figure}

We benchmark our code for FETI-DP local matrices in various mesh structures. We consider a structured and an unstructured mesh in a cube geometry geometry. The structured cube mesh uses an 8 node cube element while the unstructured cube mesh uses a 4 node tetrahedral element to discretize the elasticity equation Eq.~\eqref{eq:NavierCauchy}. Figures~\ref{fig:structuredMesh} and~\ref{fig:unstructuredMesh} show a sample mesh for the structured and the unstructured cubes respectively.

\begin{figure}[htbp]
	\centering
		\includegraphics[width=200pt,height=200pt]{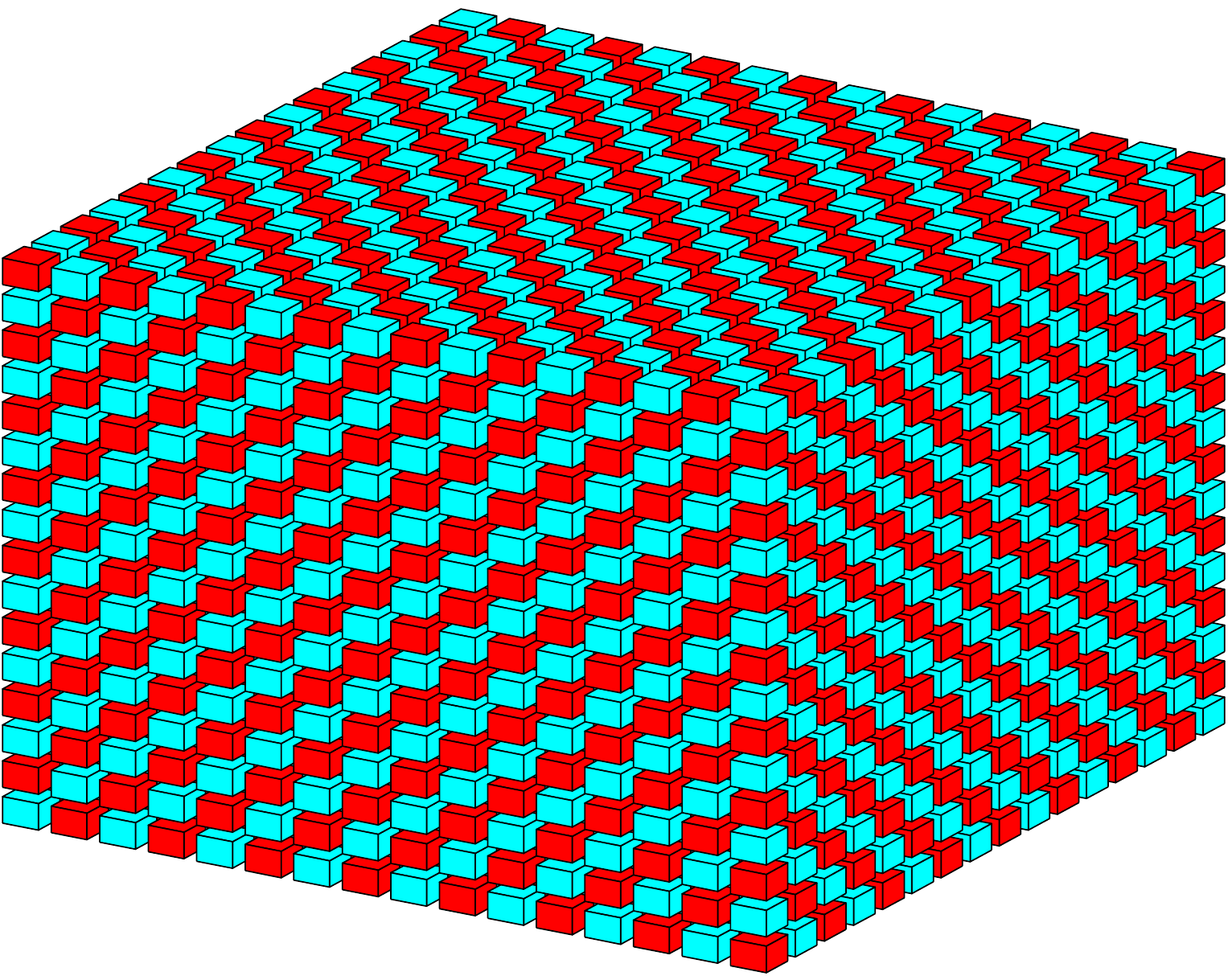}
	\caption{A sample subdomain configuration in a cube geometry. Each colored block represent a subdomain.}
	\label{fig:subdomains}
\end{figure}

We also apply our solver to FETI-DP coarse problem matrices arising from solving the elasticity equation in a unit cube geometry. Factorization of the coarse matrix for problems where the coarse matrix is large is expensive and might become the bottleneck of the FETI-DP solver. As a result, we are interested in accelerating the solve process and decreasing the memory footprint of factorizing such matrices. Figure~\ref{fig:subdomains} show a typical subdomain configuration in the unit cube.

\section{Numerical Results}

In this section we show numerical benchmarks for  the matrices described in Sections ~\ref{sec:cylinderResults} and~\ref{sec:FETIResults} respectively.  As described in Section~\ref{sec:directIterative}, we use the accelerated multifrontal solver at low accuracies as a preconditioner to the GMRES iterative method. We compare this approach to conventional preconditioners, namely the diagonal and the incomplete LU (ILU) preconditioner as well as the conventional multifrontal algorithm.

We implemented our code in C\verb|++| and used the Eigen C\verb|++| library for linear algebra operations. The incomplete LU algorithm is the incomplete LU with the dual-thresholding implementation from the SPARSEKIT package~\cite{ILUCode}.

\subsection{Elasticity Problem for a Cylinder Head Geometry}

\begin{figure}[htbp]
\centering
\subfigure[Convergence Analysis]{
\includegraphics{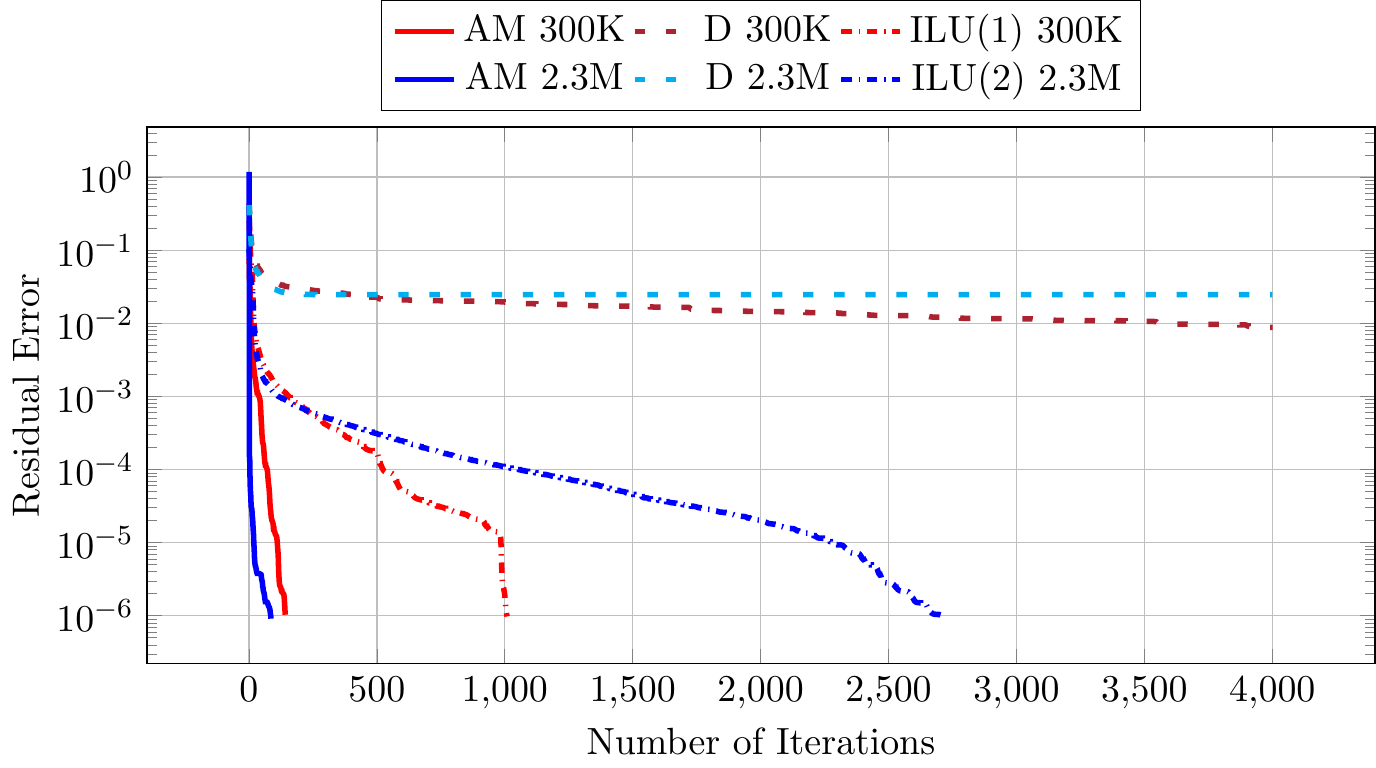}
\label{fig:cHeadIter}
}
\subfigure[RunTime]{
\includegraphics{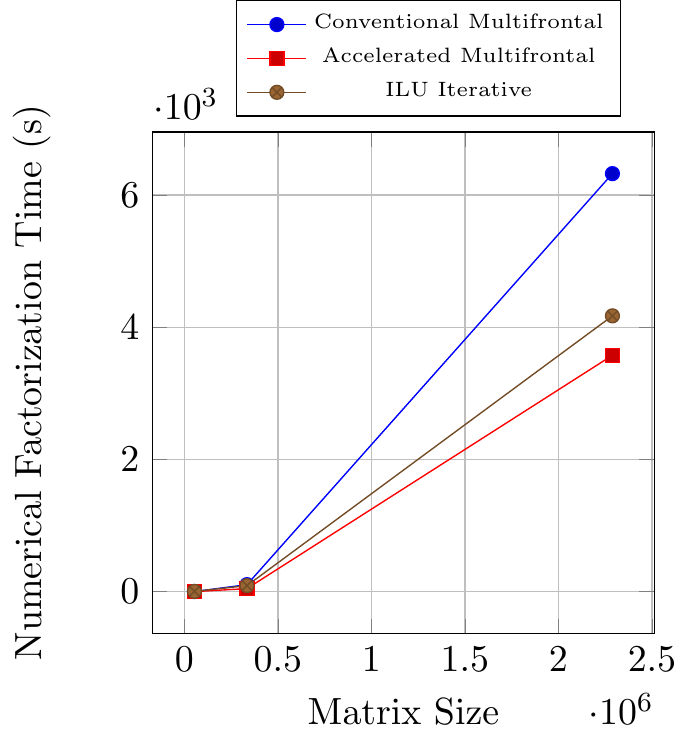}
\label{fig:cHeadTime}
}
\subfigure[Memory Consumption]{
\includegraphics{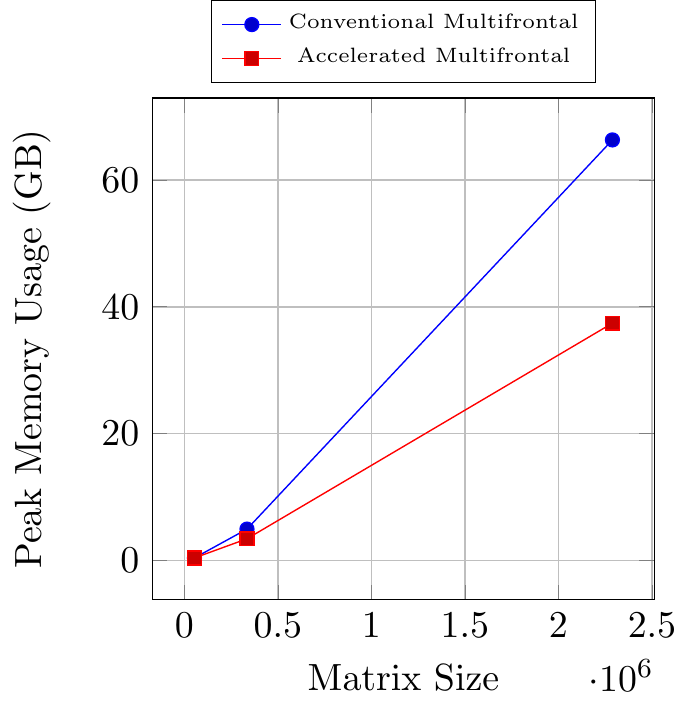}
\label{fig:cHeadMem}
}

\caption{Convergence, runtime and memory consumption analysis for the unstructured cylinder head mesh. AM stands for accelerated multifrontal preconditioner and D stands for the diagonal preconditioner. For detailed code parameters see Table~\ref{table:summaryNum}.}
\label{fig:cylinderResults}
\end{figure}

Figure~\ref{fig:cHeadIter} shows the convergence of the accelerated multifrontal preconditioner against the conventional diagonal and ILU preconditioners for the cylinder head geometry. As can be seen, since the problem is relatively difficult, the diagonal preconditioners fails to converge and one needs to work with the parameters of the ILU preconditioner in order to achieve convergence. Furthermore, the accelerated multifrontal preconditioner converges much faster compared to both diagonal and ILU preconditioners. 

Figures~\ref{fig:cHeadTime} and~\ref{fig:cHeadMem} show the run time and memory consumption comparison between the conventional multifrontal, the accelerated multifrontal and the incomplete LU iterative schemes respectively. The accelerated multifrontal algorithm has a lower runtime compared to the conventional multifrontal and the ILU algorithm. 

\subsection{FETI-DP Solver for a 3D Elasticity Problem}

\subsubsection{FETI-DP Local Problems}
\begin{figure}[htbp]
\centering
\subfigure[Convergence Analysis]{
\includegraphics{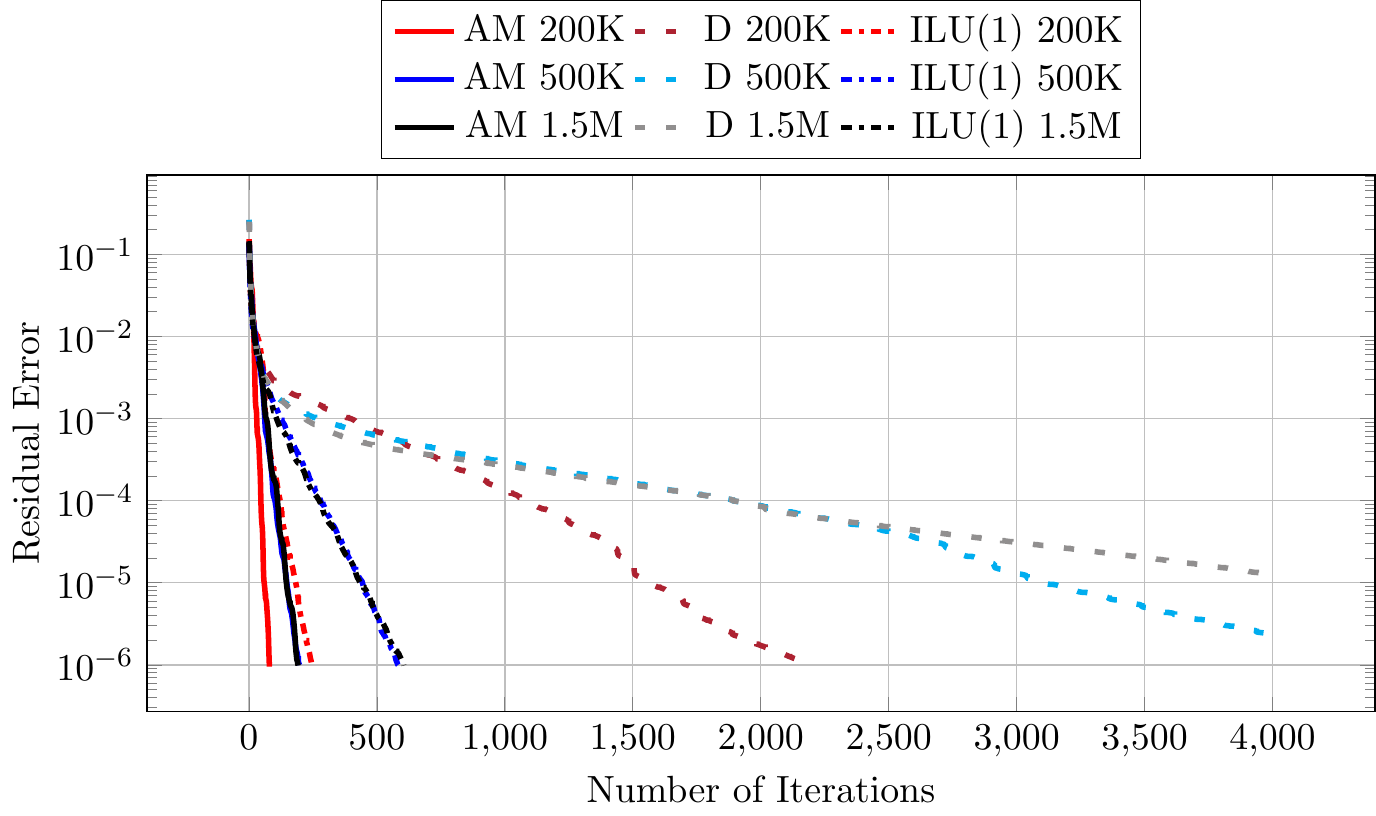}

\label{fig:structCubeIter}
}
\subfigure[Run Time]{
\includegraphics{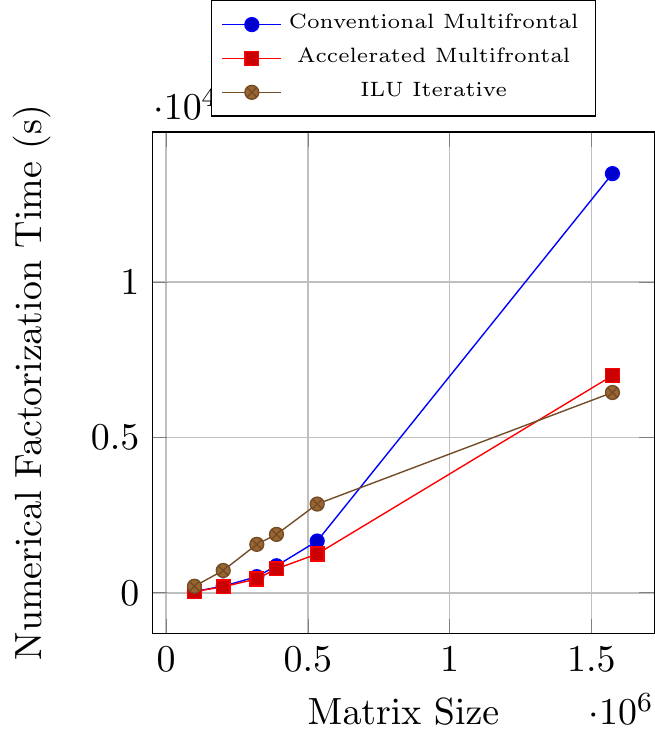}
\label{fig:structCubeTime}
}
\subfigure[Memory Consumption]{
\includegraphics{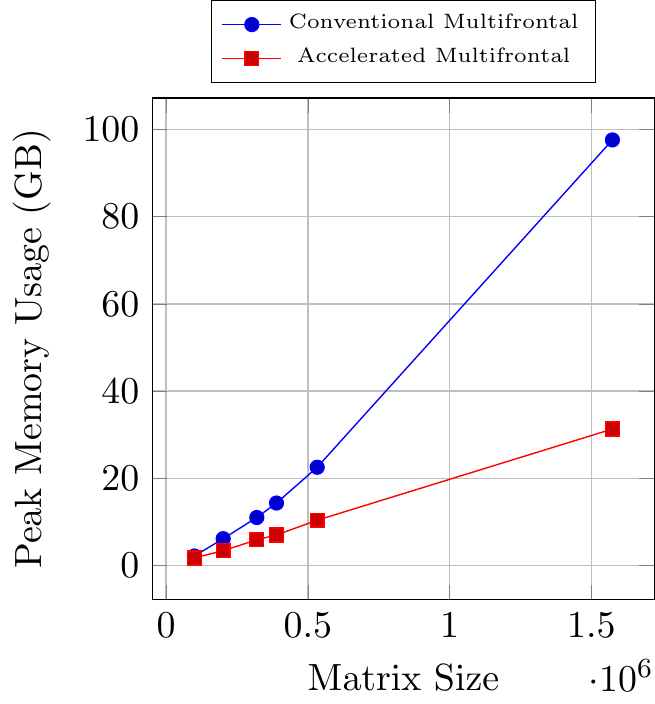}
\label{fig:structCubeMem}
}
\caption{Convergence, runtime and memory consumption analysis for FETI-DP local matrices arising from the structured cube mesh. AM stands for accelerated multifrontal preconditioner and D stands for the diagonal preconditioner. For detailed code parameters see Table~\ref{table:summaryNum}.}
\label{fig:strcuturedResults}
\end{figure}

Figure~\ref{fig:structCubeIter} compares the convergence of the accelerated multifrontal method with traditional preconditioners for the structured cube mesh local problem. As this problem is a relatively simple problem, both the diagonal and ILU preconditioners converge without too many iterations. However, the accelerated multifrontal method still has the highest convergence rate amongst the benchmarked algorithm. Figures~\ref{fig:structCubeTime} and~\ref{fig:structCubeMem} show the runtime and memory consumption comparison for the structured cube local problem. The fast multifrontal algorithm has a significantly lower runtime and memory consumption compared to the conventional multifrontal algorithm. However, because of the relative simplicity of this problem, the ILU algorithm is competitive in terms of factorization time.

\begin{figure}[htbp]
\centering
\subfigure[Convergence Analysis]{
\includegraphics{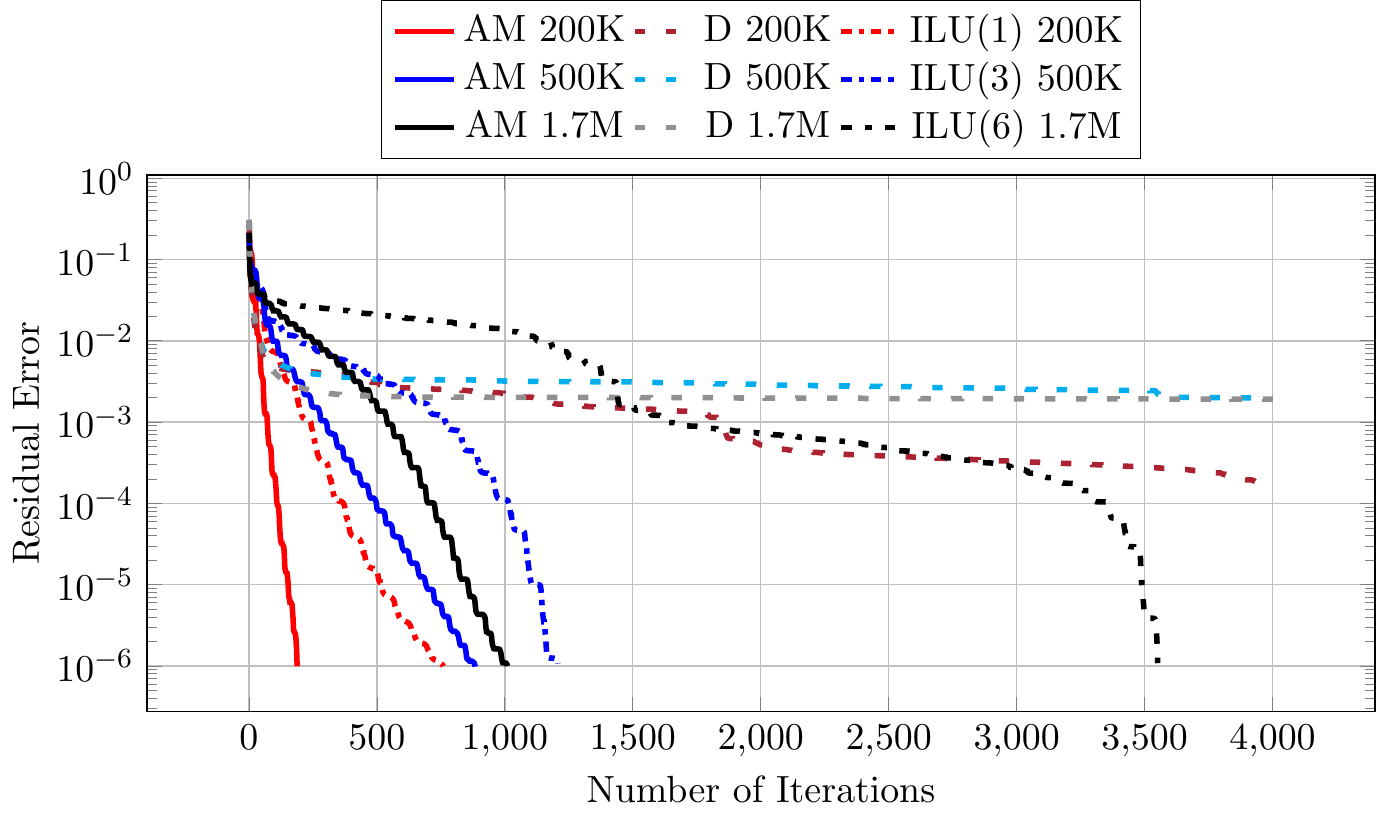}

\label{fig:unstructCubeIter}
}
\subfigure[Run Time]{
\includegraphics{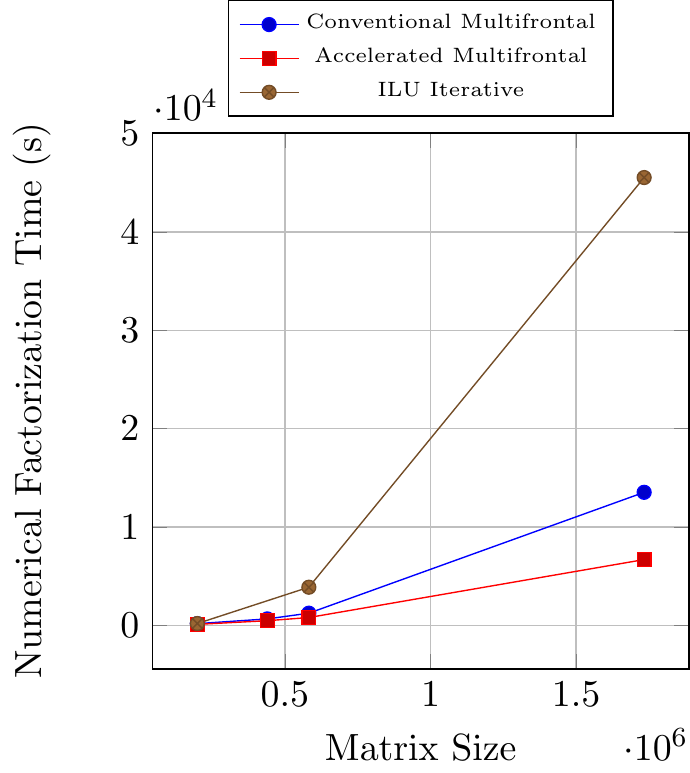}

\label{fig:unstructCubeTime}
}
\subfigure[Memory Consumption]{
\includegraphics{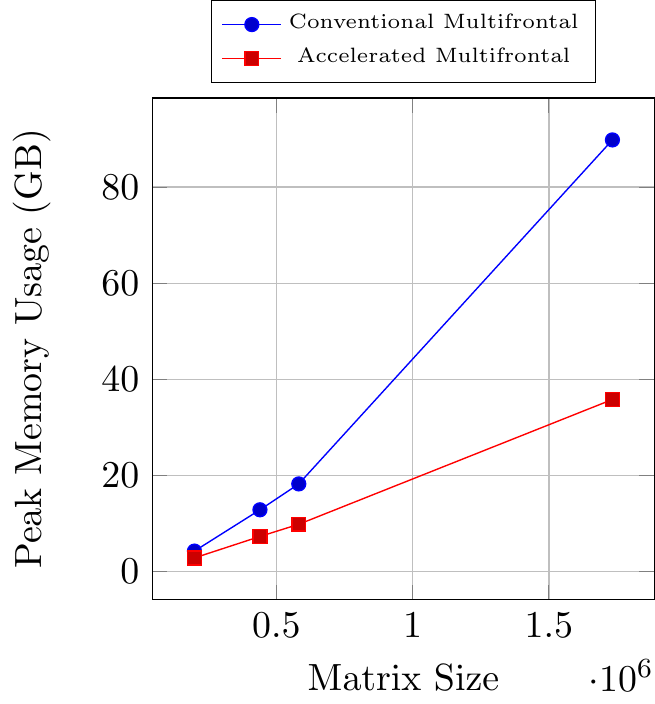}

\label{fig:unstructCubeMem}
}
\caption{Convergence, runtime and memory consumption analysis for FETI-DP local matrices arising from the unstructured cube mesh. AM stands for accelerated multifrontal preconditioner and D stands for the diagonal preconditioner. For detailed code parameters see Table~\ref{table:summaryNum}.}
\label{fig:unstructuredResults}
\end{figure}

Figure~\ref{fig:unstructCubeIter} shows the convergence rate of the accelerated multifrontal, diagonal and the ILU preconditioner for the unstructured cube mesh local problem. As can be seen, this problem is the most complicated and difficult. Not only does the diagonal preconditioner fail to converge, but the input parameters of ILU need to be increased significantly in order to achieve convergence. Figure~\ref{fig:unstructCubeTime} shows that the ILU preconditioner is significantly slower than both the multifrontal and accelerated multifrontal solver. The accelerated multifrontal solver is the fastest among all conventional algorithms, and Figure~\ref{fig:unstructCubeMem} shows that it also reduces the memory requirements. 

\subsubsection{FETI-DP Coarse Problems}
\begin{figure}[htbp]
\centering
\subfigure[Convergence Analysis]{
\includegraphics{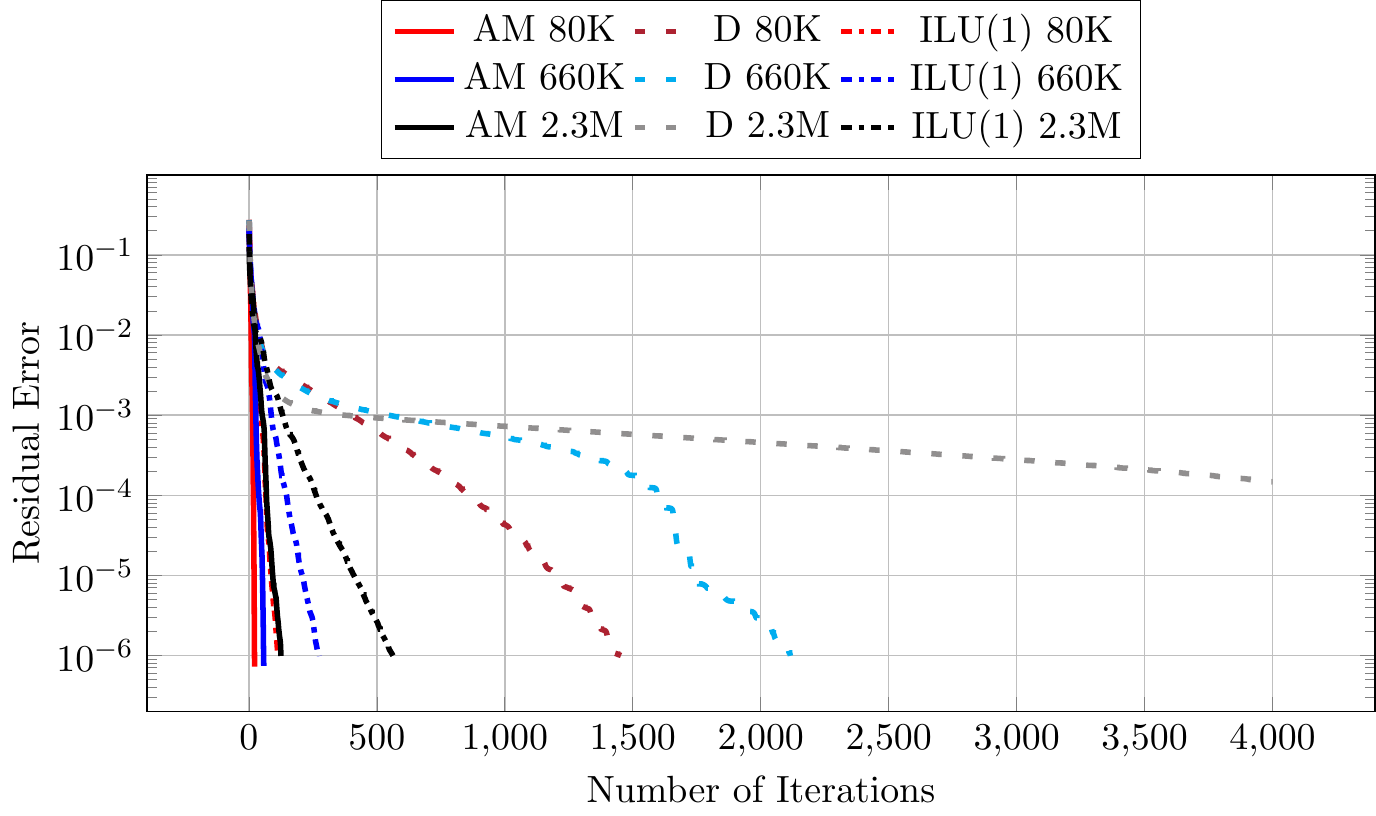}
\label{fig:elasticitySIter}
}
\subfigure[RunTime]{
\includegraphics{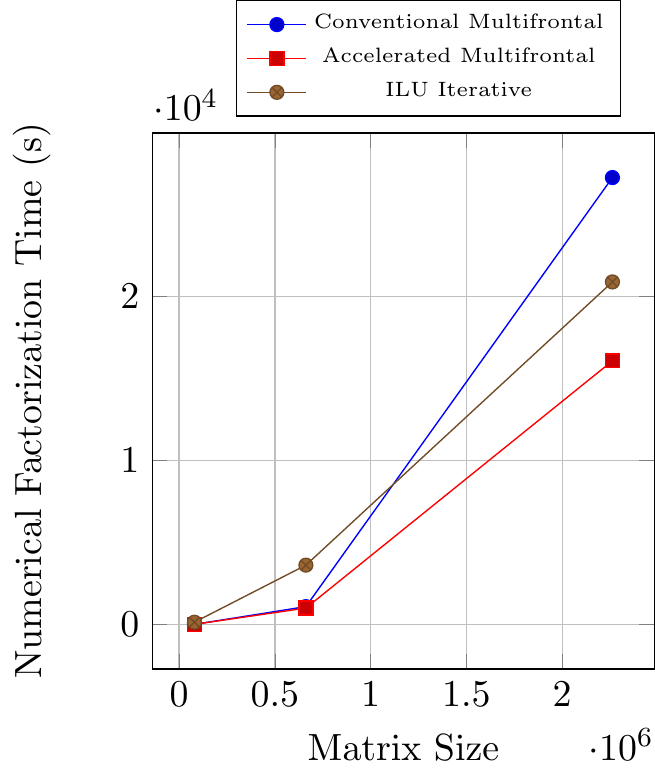}

\label{fig:elasticitySTime}
}
\subfigure[Memory Consumption]{
\includegraphics{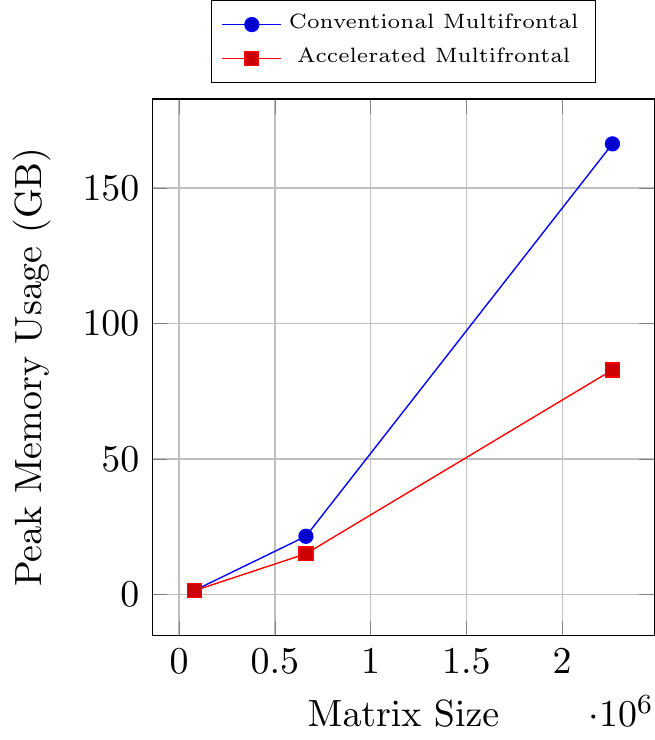}

\label{fig:elasticitySMem}
}

\caption{Convergence, runtime and memory consumption analysis for FETI-DP coarse matrices arising from the discretization of the elasticity equation in a structured cube mesh. AM stands for accelerated multifrontal preconditioner and D stands for the diagonal preconditioner. For detailed code parameters see Table~\ref{table:summaryNum}. The benchmark matrices correspond to dividing the unit cube into $16^3$, $32^3$, and $48^3$ subdomains. The size of each subdomain is $8\times8\times8$ elements. The coarse matrix is based on the displacement of corners of subdomains and the average augmentation for displacements and rotations on the faces.}
\label{fig:elasticitySResults}
\end{figure}

Figure~\ref{fig:elasticitySIter} shows the convergence rate of the accelerated multifrontal method and the conventional preconditioners for a coarse FETI-DP problem in a cube geometry. Figure~\ref{fig:elasticitySTime} shows that the accelerated multifrontal method is faster than both the conventional multifrontal and the ILU algorithms. Furthermore, as can be seen in  Figure~\ref{fig:elasticitySMem}, the memory consumption of the accelerated multifrontal method is significantly lower compared to the conventional algorithm.

\begin{figure}[htbp]
\centering
\subfigure[Convergence Analysis]{
\includegraphics{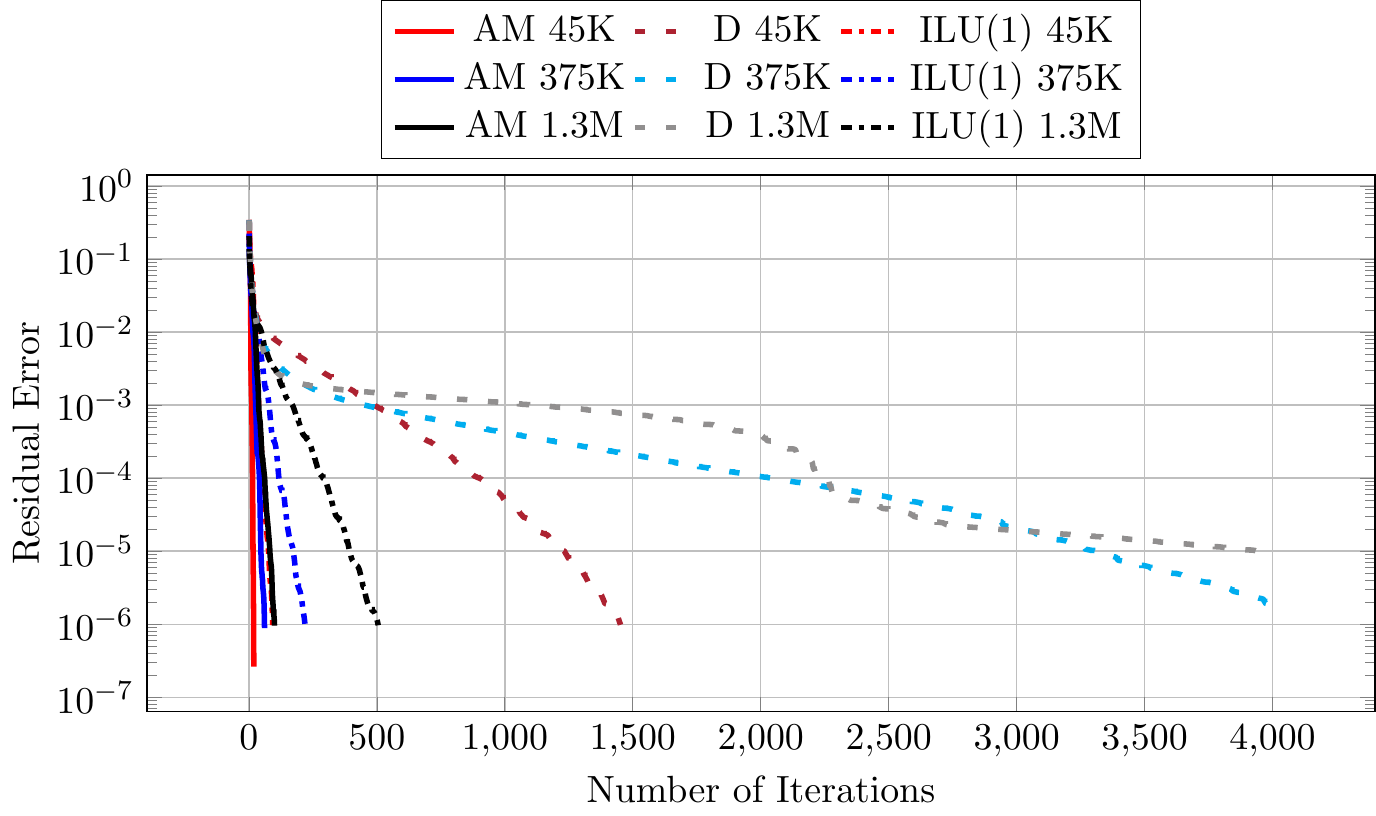}

\label{fig:elasticityStIter}
}
\subfigure[RunTime]{
\includegraphics{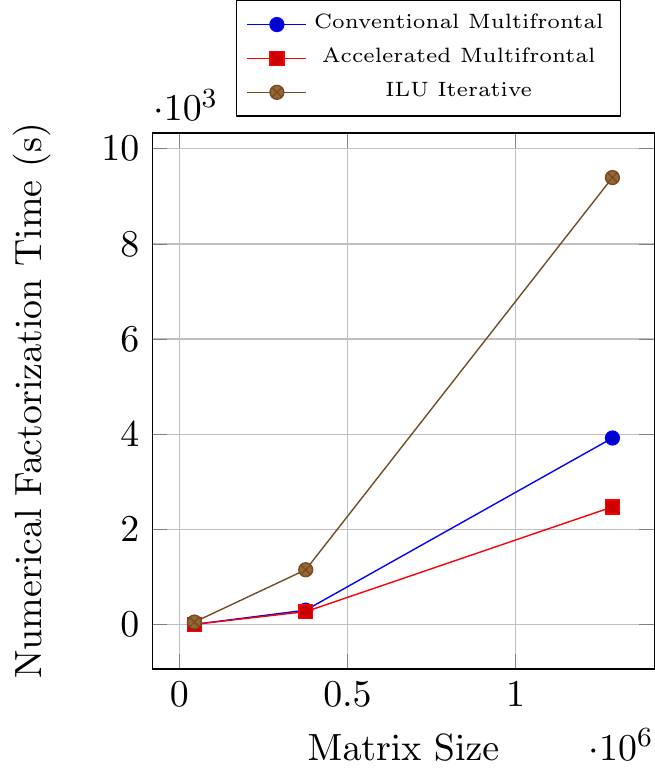}

\label{fig:elasticityStTime}
}
\subfigure[Memory Consumption]{
\includegraphics{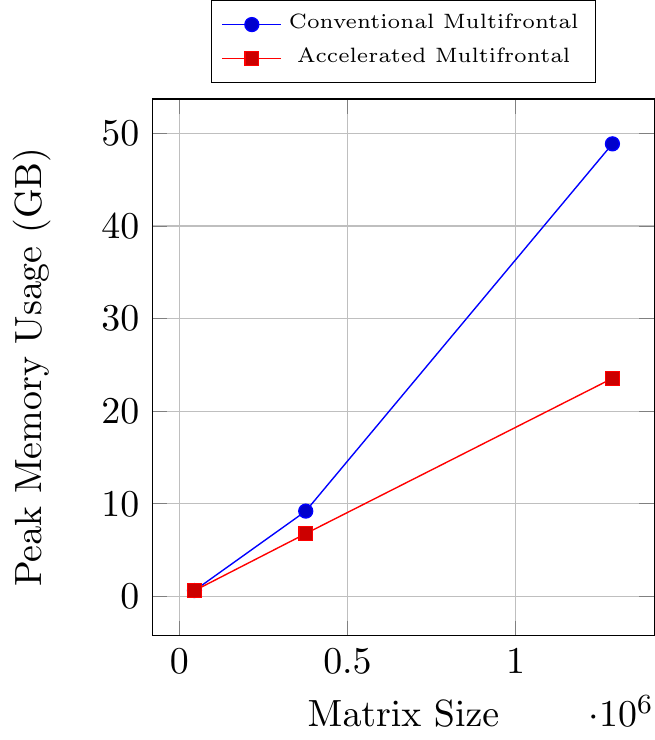}

\label{fig:elasticityStMem}
}

\caption{Convergence, runtime and memory consumption analysis for FETI-DP coarse matrices arising from the discretization of the Elasticity equation in a structured cube mesh. AM stands for accelerated multifrontal preconditioner and D stands for the diagonal preconditioner. For detailed code parameters see Table~\ref{table:summaryNum}. The benchmark matrices correspond to dividing the unit cube into $16^3$, $32^3$ and $48^3$ subdomains. The size of each subdomain is $8\times8\times8$ elements. The coarse matrix is based on corners of subdomains and the average augmentation for displacements without rotations on the faces.}
\label{fig:elasticityStResults}
\end{figure}

Figure~\ref{fig:elasticityStIter} compares the convergence rates of the accelerated multifrontal method with ILU and diagonal preconditioning schemes for a coarse FETI-DP problem that only includes translational degrees of freedoms at the corner of the subdomains. As can be seen in Figure~\ref{fig:elasticityStTime}, the ILU algorithm is significantly slower compared to both the conventional and the accelerated multifrontal schemes. Figures~\ref{fig:elasticityStTime} and~\ref{fig:elasticityStMem} show that not only the accelerated multifrontal method is faster, but also it consumes much less memory compared to the conventional multifrontal algorithm.

\subsection{Summary}
\label{sec:summary}

Table~\ref{table:summaryNum} shows a detailed summary of all the benchmark cases. As can be seen in all cases, the GMRES solver with the accelerated multifrontal preconditioner converges after few iterations which shows the effectiveness of the developed algorithm. Furthermore, in almost all cases, we observe a speedup and memory saving of up to 3x compared to the conventional multifrontal solver. We were not able to benchmark larger cases due to memory limitations. However, one can observe that both the speedup and memory saving become more significant as the matrix size grows. That is, a very high speedup and memory saving can be achieved for very large matrices.

\begin{figure}[htbp]
\centering
\subfigure[Number of Iterations for Structured Meshes]{
\includegraphics{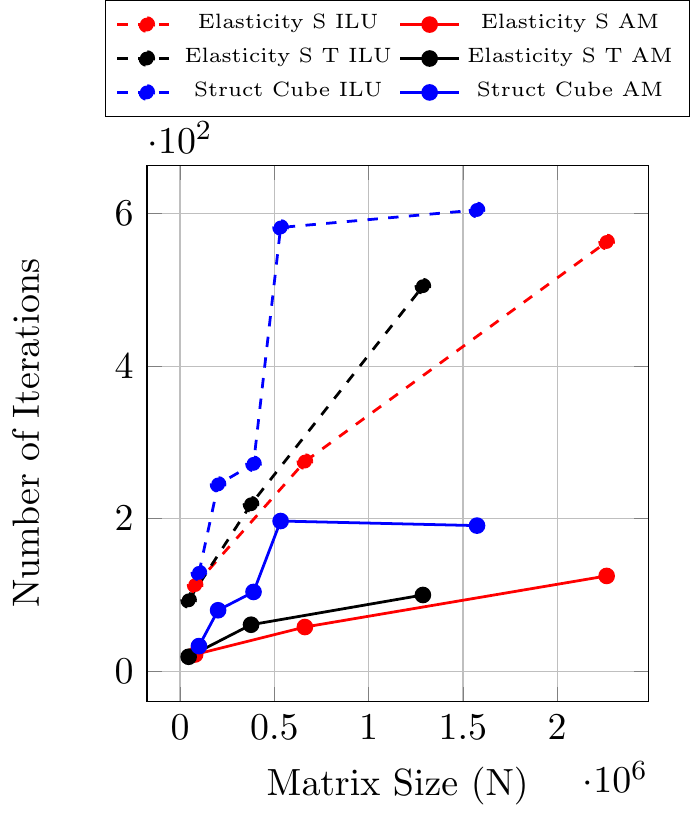}

\label{fig:numIterStruct}
}
\subfigure[Number of Iterations for Unstructured Meshes]{
\includegraphics{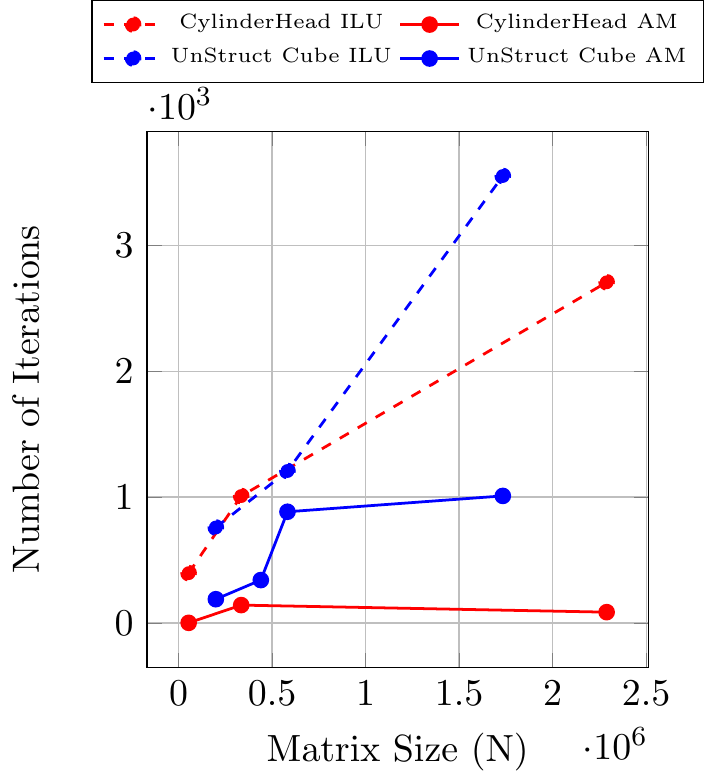}

\label{fig:numIterUnstruct}
}

\subfigure[Convergence Analysis]{
\includegraphics{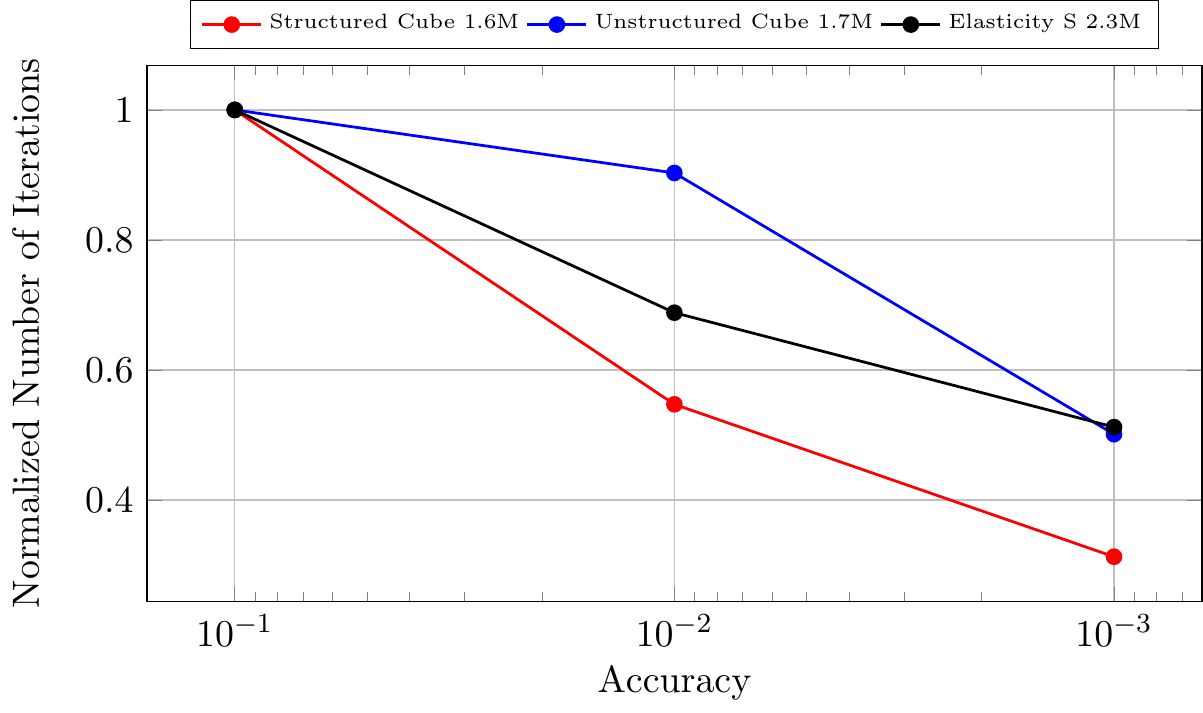}

\label{fig:numIterAcc}
}
\caption{Number of iterations vs.\ matrix size and solver accuracy for a variety of problems. a) Number of iterations vs.\ matrix size for problems with a structured cube mesh. b) Number of iterations vs.\ matrix size for problems with an unstructured mesh. c) Normalized number of iterations vs.\ fast solver accuracy. Number of iterations has been normalized by the number of iterations at the accuracy of $10^{-1}$.}
\end{figure}

Figures~\ref{fig:numIterStruct} and~\ref{fig:numIterUnstruct} compare the number of iterations for the ILU and the accelerated multifrontal preconditioner for the benchmarked structured and unstructured meshes respectively. As both figures show, not only the accelerated multifrontal preconditioner has less number of iterations compared to ILU, but also the number of iterations does not grow significantly with matrix size. This is another important advantage of the accelerated multifrontal preconditioner that makes it more favorable for parallel implementations compared to ILU. This is because fewer number of iterations results in fewer matrix vector products which requires less communication between the nodes which ultimately results in a higher speedup.

Figure~\ref{fig:numIterAcc} shows that one can significantly decrease the number of iterations required by the accelerated multifrontal preconditioner by simply  decreasing the accuracy parameter (and in turn, increasing the depth parameter). However, this results in an increase in the off-diagonal rank and the decrease of speedup and memory savings. This shows that one can fine-tune the code parameters based on their available resources and their desired convergence rates.

\begin{table}[htbp]
\centering
	\scalebox{0.7}{
	\begin{tabular}{|c||c||c|c|c|c|c|c|c|c|c|c|c|c|c|c|}
	\hline
	 \multirow{3}{*}{Matrix} & \multirow{3}{*}{Mesh} & \multirow{3}{*}{Matrix}& \multicolumn{2}{c}{Conventional}\vline & \multicolumn{6}{c}{Accelerated} \vline&\multicolumn{3}{c}{GMRES}\vline& \multirow{3}{*}{Speed} & \multirow{3}{*}{Mem} \\ \cline{12-14}
	  & & & \multicolumn{2}{c}{Multifrontal}\vline & \multicolumn{6}{c}{Multifrontal} \vline&{D}&\multicolumn{2}{c}{ILU}\vline& \multirow{3}{*}{-up}&  \\  \cline{4-14} 
	  Type& Type&{Size} & Fact & Mem& Fact & Mem& Num&\multicolumn{3}{c}{Parameters}\vline&Num& Num&\multirow{2}{*}{k}&&Saving\\ \cline{9-11}
	  &  &  & (s) & (GB)& (s) & (GB) &Iter&$n_c$&$\epsilon$&$d$ &Iter&Iter&&&\\ 	\hline
	  \multirow{2}{*}{Stiffness}&\multirow{2}{*}{C Head}
	 &330K  &1.08e2&4.92&4.78e1&3.39  &142&3K&1e-1&1 &x&1009&1&\bf2.26&1.645\\
	 &&2.30M*&6.33e3&66.34&3.58e3&37.42  &86&10K&1e-2&5&x&2709&2&\bf1.77&\bf1.77\\ \cline{1-16}
	  \multirow{9}{*}{FETI }&\multirow{5}{*}{Str}
	   &100K  &4.32e1&2.16  &3.88e1&1.72&33&3K  &1e-1&1&813&129&1&1.11&1.25\\ 
	  \multirow{9}{*}{Local} &\multirow{5}{*}{Cube}&200K  &2.14e2&6.13  &1.94e2&3.38&80&3K &1e-1&1&2194&245&1&1.10&\bf1.81\\
	  &&320K  &5.20e2&10.99&4.48e2&5.87&81&3K&1e-1&1&1059&216&1&1.16&\bf1.87\\
	  &&390K  &8.69e2&14.30&7.74e2&6.99&104&3K&1e-1&1&2759&272&1&1.12&\bf2.05\\
	  &&530K  &1.67e3&22.52&1.25e3&10.36&197&3K&1e-1&1&x&582&1&1.34&\bf2.17\\
	  &&1.57M&1.35e4&97.62&6.99e3&31.29&191&3K&1e-1&1&x&605&1&\bf1.93&\bf3.11\\\cline{2-16}
	 &\multirow{3}{*}{Uns}
	 &200K  &1.74e2  &4.26  &9.93e1&2.84&189&4K&1e-1&1 &x&760&1&\bf1.75&1.47\\
	 &\multirow{3}{*}{Cube}&440K  &6.55e2  &12.85&4.73e2&7.30&341&5K&1e-1&1&x&1084&3&1.38&\bf1.76\\
	 &&580K*  &1.24e3  &18.24&7.95e2&9.84&884&6K&1e-1&1 &x&1209&3&\bf1.56&\bf1.85\\
	 &&1.73M*&1.35e4  &89.82&6.69e3&35.81&1010&10K&1e-1&1&x&3553&6&\bf2.02&\bf2.05\\ \cline{1-16}
	  \multirow{5}{*}{FETI}&\multirow{3}{*}{Elasticity}
	   &80K*  &1.90e1  &1.51  &2.14e1&1.42&22&3K&1e-1&1 &1454&113&1&0.88&1.06\\ 
	   \multirow{5}{*}{Coarse}&&660K*  &1.10e3&21.50  &1.01e3&15.04&58&18K&1e-1&1 &2116&275&1&1.08&1.43\\
	 &&2.26M*&2.72e4  &166.35&1.61e4&82.89&125&16K&1e-1&1&x&563&1&\bf1.69&\bf2.00\\ \cline{2-16}
	  &\multirow{2}{*}{Elasticity}
	   &45K*  &5.37e0  &0.63  &5.89e0&0.62&19&2K&1e-1&1 &1452&93&1&0.91&1.02\\ 
	 &\multirow{2}{*}{T}&375K*  &3.04e2&9.21  &2.72e2&6.78&61&5K&1e-1&1 &x&219&1&1.12&1.36\\
	 &&1.30M*&3.92e3  &48.87&2.47e3&23.55&100&10K&1e-1&1&x&505&1&\bf1.59&\bf2.08\\ \hline
	\end{tabular}}
  \caption{Summary of solver accuracy, speed and memory consumption for various benchmark cases. All timings are measured in seconds and memory usage is measured in Giga Bytes. {\bf FETI local:} FETI-DP local matrices. {\bf FETI coarse:} FETI-DP coarse problem matrices. T stands for a matrix where, in the coarse FETI-DP matrix, we consider only the translational degrees of freedom for the corner nodes. {\bf Stiffness:} Regular finite-element stiffness matrices. `Str' refers to structured meshes and `Uns' refers to unstructured meshes. All results are obtained using the GMRES iterative method with a termination tolerance of $10^{-6}$. Columns `AM' ,`D'  and `ILU' refer to the GMRES method with the accelerated multifrontal, diagonal and incomplete LU preconditioners respectively. The letter `x' indicates that the respective iterative method did not converge within 4,000 iterations. The parameters column reflects the code parameters that were used in obtaining the results. $n_{c}$: size threshold for converting from dense linear algebra to HODLR matrix operations. $\epsilon$: error tolerance used in the low-rank approximations. $d_{BDLR}$: depth parameter in the BDLR low-rank approximation method. `k' is a parameter that is used to identify the amount of fill-in in the ILU scheme. That is, the fill-in in each row is set to $\frac{k \cdot NNZ}{N}+1$ where $N$ and $NNZ$ are the matrix size and the number of non zeros in the matrix respectively. Speed-up denotes the speed up in numerical factorization phases, whereas memory savings denotes the overall memory savings compared to the conventional multifrontal method. The symbol `*' in front of the matrix name depicts that the matrix has been scaled such that the norm of the largest entry in each row is $1$.}	
\label{table:summaryNum}	
\end{table}

\section {Conclusion}
We have developed a black-box fast and robust linear solver for finite element matrices arising from the discretization of elliptic PDEs. As we've shown, this solver is advantageous both in terms of running time and memory consumption compared to both the conventional multifrontal direct solve algorithm and the ILU preconditioner. Furthermore, not only our solver is faster in terms of factorization time, but also it results in less number of iterations compared to ILU.

The examples presented here were run on a single core machine, and are limited in size by the amount of memory available on a single computer node. A parallel implementation would have allowed us to run much larger test cases. Since the speed-up improves with $N$, this would have allowed us to demonstrate even greater speed-ups, in particular compared to ILU. This will be, however, done in a future publication. 

This solver can be used at low-accuracy as a preconditioner, or as a standalone direct solver at high accuracy. The current scheme relies on the assumption that all off-diagonal blocks are low-rank. In practice, this implies that the rank required to reach the desired accuracy may become somewhat large, leading to a loss in efficiency. This can be remedied using more complex algorithms such as~\cite{ambikasaran2014ifmm}. Such algorithms are currently under development for the case of sparse matrices. Despite this limitation, the class of methods presented does yield significant improvements over the current state of the art.

One advantage of the method presented here is its relative simplicity. For example, by removing the requirement to form a nested low-rank basis across levels, we can simplify the implementation and algorithm significantly. This is in contrast with the HSS class of methods for example~\cite{sheng2007algorithms}. Despite this simplification, the HODLR scheme has a computational cost in $O(N^{4/3})$, whereas HSS-based schemes scale like $O(N^{4/3} \log N)$~\cite{randomizedMF}.

We finally point out that the algorithms presented here are very general and robust. They can be applied to a wide-range of problems in a black-box manner. This was demonstrated in part in this manuscript. 

\section*{Acknowledgments}
The authors would like to acknowledge Prof.\ Charbel Farhat, Dr.\ Philip Avery and Dr. Jari Toivanen for providing us with the FETI test matrices. We also want to thank Profs.\ Pierre Ramet and Mathieu Faverge for their useful input and suggestions in this project.

Part of this research was done at Stanford University, and was supported in part by the U.S.\ Army Research Laboratory, through the Army High Performance Computing Research Center, Cooperative Agreement W911NF-07-0027. This material is also based upon work supported by the Department of Energy under Award Number DE-NA0002373-1.

\FloatBarrier

\FloatBarrier

\bibliographystyle{elsarticle-harv}

\bibliography{sparseSolver}

\end{document}